\documentclass[a4paper,11pt,fleqn]{article}
\usepackage[utf8]{inputenc}
\usepackage[T1]{fontenc}
\usepackage[english]{babel}
\usepackage{soul}
\usepackage[pdftex]{graphicx}
\usepackage{geometry}
\usepackage{tabularx}
\usepackage{multirow}
\usepackage{float}
\usepackage{charter}

\usepackage{titlesec}
\titleformat*{\section}{\Large\bfseries\color{black}}
\titleformat*{\subsection}{\large\itshape\color{black}}
\titleformat*{\subsubsection}{\small\itshape\color{black}}

\usepackage{tikz}
\usetikzlibrary{positioning}
\usetikzlibrary{calc,through,backgrounds}
\usepackage{pgfplots}
\pgfplotsset{compat=newest}
\usepackage{changepage}

\usepackage{algorithm,,algpseudocode}

\usepackage{hyperref}
\usepackage{authblk}

\usepackage{amsfonts}
\usepackage{amsmath}
\usepackage{amsthm}
\usepackage{amsmath}
\usepackage{amssymb}
\usepackage{mathrsfs}
\usepackage{booktabs}
\usepackage{siunitx}
\usepackage{thmtools}
\usepackage{ulem}

\usepackage{xcolor}
\usepackage{times}
\usepackage{tikz}
\usepackage{lipsum}

\declaretheoremstyle[
    bodyfont=\normalfont\color{red},
    headfont=\color{red}
]{styleattention}

\declaretheoremstyle[
    spacebelow=1em
]{styleremarque}

\declaretheoremstyle[
    spaceabove=-6pt,
    spacebelow=6pt,
    headfont=\normalfont\bfseries,
    bodyfont=\normalfont,
    postheadspace=1em,
    qed=$\Box$,
]{mystyle}

\declaretheorem[thmbox=M,numberwithin=section,title=Definition]{definition}
\declaretheorem[thmbox=M,sibling=definition]{proposition}

\declaretheorem[thmbox=M,sibling=definition,title=Lemma]{lemme}

\declaretheorem[style=styleremarque,sibling=definition,title=Remark]{remarque}

\declaretheorem[name={}, style=mystyle, unnumbered, title=Proof]{preuve}

\let\d\undefined
\newcommand{\d}{\; \mathrm{d}}

\newcommand{\1}{{\mathchoice {\rm 1\mskip-4mu l} {\rm 1\mskip-4mu l}{\rm 1\mskip-4.5mu l} {\rm 1\mskip-5mu l}}}

\newcommand{\N}{\mathbb{N}}

\newcommand{\prox}{\text{prox}}
\let\xp\undefined
\newcommand{\xp}{\xi_p} 
\newcommand{\bp}{\zeta_p} 
\newcommand{\zp}{\omega_p} 
\newcommand{\X}{\mathscr{X}} % Estimation space
\newcommand{\Y}{\mathscr{Y}} % Observation space
\newcommand{\yd}{y^\delta} % Observation space

\title{Inversion of Integral Models: a Neural Network Approach}
\author{E. Chouzenoux, C. Della Valle, and J.-C. Pesquet \footnote{
C. Della Valle (corresponding author) is with Université de Paris, Paris Sorbonne Université.\\
E. Chouzenoux and J.-C. Pesquet are with Universit\'e Paris-Saclay, Inria, CentraleSupélec, Center for Visual Computing.
The work by J.-C. Pesquet was supported by Institut Universitaire de France and the ANR Chair in Artificial Intelligence BRIDGEABLE.
}}
\date{}

\begin{document}

\maketitle
%%%%%%%%%%%%%%%%%%%%%%%%%%%%%%%%%%%%%%%%%%%%%%%%%%%%%%%%%%%%%%%%%%%%%%%%%%%%%%%%%%%%%%%%%%%%%%%%%%%%%%%%%%%%%%%%%%%%%%%%%%%%%%%%%
\begin{abstract}
    We introduce a neural network architecture to solve inverse problems linked to a one-dimensional integral operator. This architecture is built by unfolding a forward-backward algorithm derived from the minimization of an objective function 
    which consists of the sum of a data-fidelity function and a Tikhonov-type regularization function. 
    The robustness of this inversion method with respect to a perturbation of the input is theoretically analyzed. 
    Ensuring robustness is consistent with inverse problem theory since it guarantees both the continuity of the inversion method and its insensitivity to small noise. The latter is a critical property as deep neural networks have been shown to be vulnerable to adversarial perturbations.
    One of the main novelties of our work is to show that the proposed network is also robust to perturbations of its bias.
    %which in our case is an input of the recovery algorithm. 
    In our architecture, the bias accounts for the observed data in the inverse problem. 
    We apply our method to the inversion of Abel integral operators, which define a fractional integration involved in wide range of physical processes. 
    The neural network is numerically implemented and tested to illustrate the efficiency of the method. 
    Lipschitz constants after training are computed to measure the robustness of the neural networks.
\end{abstract}

%%%%%%%%%%%%%%%%%%%%%%%%%%%%%%%%%%%%%%%%%%%%%%%%%%%%%%%%%%%%%%%%%%%%%%%%%%%%%%%%%%%%%%%%%%%%%%%%%%%%%%%%%%%%%%%%%%%%%%%%%%%%%%%%%
\section{Introduction}
\paragraph{Inverse problem}
In this article, we are interested in 1D operators of the form
\begin{equation}
	   \label{def:T}
  	\begin{array}{cc|ccc}
  	 T & : & \X  & \to & \Y \\
  	   &   & x   & \to & 
	 \displaystyle y(t) = \int_0^1 k(t,s) x(s) \d s \; .\\
  	\end{array}
 \end{equation}
 Hereabove, $\X$ and $\Y$ are functional Hilbert spaces, typically $\X=\Y = L^2(0,1)$, 
$k\in L^2([0,1]^2)$ and $T$ is a linear compact operator. 
A large variety of inverse problems consist of inverting convolution operators such as signal/image restoration~\cite{Bertero2009,Chouzenoux2012}, tomography~\cite{Wu2020}, Fredholm equation of the first kind~\cite{Arsenault2017} or inverse Laplace transform~\cite{Cherni2017}.
In this work,
we focus on the inversion of the Abel integral, for which
the kernel is of the form 
 \[
  k(t,s) = \ell (t,s) (t-s)^{a-1} \delta_{s\, \leq \, t }\; ,
 \]
 where $a$ is a real positive number,
 $\ell$ is a continuous function, differentiable and decreasing on its second variable,
 and $\delta_{s\, \leq \, t }$ is equal to one if $s\, \leq \, t $ and zero otherwise.
% %
The inverse problem investigated in this article is the following:
given $y \in \Y$, we seek for $x \in \X$ such that
 \[
 y = Tx \; .
 \]
In addition, we consider the case when the data $y$ are corrupted with measurement errors or noise.
We model an additive noise as follows: we call the upper
  bound for the noise level $\delta>0$ measured in $\Y^\delta$ with $\Y \subset \Y^\delta$.
  Here, $\Y^\delta$ is the Hilbert space $H^{-s}(0,1)$ defined in~\cite{DINEZZA2012521}. Typically, $s=0$ for a deterministic noise and $s=1/2$ for a deterministic equivalent of a Gaussian white noise~\cite{nussbaum1996asymptotic}.
Solving the inverse problem in the presence of noisy data then amounts to finding $x^\delta \in  \Y$ such that
\begin{equation}
    \label{def:InversePb}
     \yd = T x^\delta \; , \quad  \text{with} \quad  \|y - \yd \|_{\Y^\delta} \leq \delta \;   .
\end{equation}
The above problem is often ill-posed i.e., 
 a solution might not exist, 
 might not be unique, 
 or might not depend continuously on the data.

\paragraph{Variational problem} 
The well-posedness of the inverse problem defined by~\eqref{def:InversePb} is retrieved by regularization.
Here we consider Tikhonov type regularization.
Let $\tau \in ]0,+\infty[$ be the regularization parameter.
Solving the inverse problem ~\eqref{def:InversePb} with such regularization, leads to the resolution of the following optimization problem
\begin{equation}
\label{def:varJ0}
\underset{x \in C}{\text{minimize}}\, 
J_{\tau}(x)
\;,
\end{equation}
where 
\begin{equation}
(\forall x \in \X)\quad  J_\tau(x) = \frac{1}{2} \| Tx - \yd \|^2 + \frac{\tau}{2} \| D^r x \|^2,
\end{equation}
$C$ is a nonempty closed convex subset  of $\X$,
and $D^{r}$ acts as a derivative operator with order $r \geq 0 $.
Often, we have an a priori of smoothness on the solution, in $H^q(0,1)$ with $q\geq 0$, 
which justifies the use of such a derivative-based regularization.
Problem \eqref{def:varJ0} is an instance of the more general problem stated above, encountered in many signal/image processing tasks:
\begin{equation}
\label{def:varJ}
\underset{x \in \X}{\text{minimize}}\,
J_{\tau}(x)+\mu\, g(x)
\;,
\end{equation}
where $\mu \in [0,+\infty[$ is an additional regularization constant and
$g$ is a proper lower-semicontinuous convex function from some Hilbert space 
$\X$ to $]-\infty,+\infty]$. Indeed, Problem \eqref{def:varJ0} corresponds to the case when $g$ is the indicator function $\iota_C$
of set $C$.

\paragraph{Neural network}
We focus our attention on seeking for a solution to the addressed inverse problem through nonlinear approximation techniques
making use of neural networks. 
Thus, instead of considering the solution to the regularized problem~\eqref{def:varJ},
we define the solution to the inverse problem~\eqref{def:InversePb}
as the output of a neural network, whose structure is similar to  a recurrent network~\cite{Yu2019}. 

Namely, by setting an initial value $x_0$, we are interested in the following $m$-layers neural network
where $m\in \mathbb{N}\setminus\{0\}$:
\begin{equation}
	\label{def:modelNN}
	\begin{cases}
	\textbf{Initialization:} \\
	\quad b_0 = T^*\yd ,\\
	\textbf{Layer $n\in \{1,\ldots,m\}$:} \\
      \quad x_n 
          = R_n(W_n x_{n-1} + V_n b_0)\;,
    \end{cases}		  
\end{equation}
where, for every $n\in \{1,\ldots,m\}$,
\begin{align}
&R_n = \text{prox}_{\lambda_n \mu_n g}\\
&W_n =  \1 - \lambda_n T^*T - \lambda_n \tau_n D^*D \\
&V_n = \lambda_n \1.
\end{align}
Hereabove, $\prox_{\varphi}$ states for the proximity operator of a lower-semicontinuous proper convex function $\varphi$ \cite[Chapter 9]{bauschke2011},
$\1$ denotes the identity operator, and for every $n\in \{1,\ldots,m\}$, $\lambda_{n}$, $\mu_{n}$,
and $\tau_{n}$ are positive constants, which are learned during training. Throughout this paper, $L^*$ denotes the adjoint of a bounded linear operator~$L$ defined
on Hilbert spaces.

Model~\eqref{def:modelNN} can be viewed as unrolling $m$ iterations of an optimization algorithm, so leading to Algorithm~\ref{algo:NN}. Note that, when $\mu_n \equiv \mu$ and $\tau_n \equiv \tau$,
 we recognize a forward-backward algorithm~\cite{combettes2005signal,combettes2011proximal} 
 applied to the variational problem~\eqref{def:varJ}.
 \begin{algorithm}
 	\caption{Proximal forward-backward splitting method} 
  \label{algo:NN}
 	\begin{algorithmic}[1]
		\State Set $x_0$ ,
 		\For {$n=1,2,\ldots,m$} 
 		    \State Set $ \lambda_n, \tau_n, \mu_n $ ,
 			\State $x_{n} = \text{prox}_{\lambda_n \mu_n g} \; \left( x_{n-1} - \lambda_n \nabla J_{\tau_n} (x_{n-1}) \right) $ ,
 		\EndFor
 		\State \textbf{return} $x_{m}$
 	\end{algorithmic} 
 \end{algorithm}

From a theoretical standpoint, there is no guarantee that such a model constitutes a regularizing family, and there is no equivalence between the regularized inverse problem~\eqref{def:varJ} and the output of Model~\eqref{def:modelNN}, since the number of iterations $m$ is fixed in advance. 
However, we can quantify the robustness of Model~\eqref{def:modelNN} to perturbations on its initialization~$x_0$ and on its bias $b_0$, by an accurate estimation of its Lipschitz constant.

\paragraph{Related works and contributions}
There has been a plethora of techniques developed to invert integrals of the form \eqref{def:T}.
Among these methods, Tikhonov-type methods are attractive from a theoretical viewpoint, especially because they provide good convergence rate as the noise level decreases, as shown in~\cite{hegland1995} or~\cite{hofmann2005}. 
However, limitations of such methods may be encountered in their implementation. 
Indeed, certain parameters such as gradient descent steps or the regularization coefficient
need to be set, as discussed in~\cite{aakesson2008} or~\cite{daun2006} for the Abel integral operator. 
The latter parameter depends on the noise level, as shown in~\cite{Engl1996}, which is not always easy to estimate.
In practical cases, methods such as the L-curve method, see~\cite{hansen1999}, can be implemented to set the regularization parameter, but they require a large number of resolutions and therefore a significant computational cost.
Moreover, incorporating constraints on the solution may be difficult in such approaches, and often reduces to projecting the resulting solution onto the desired set.
These points justify the use of a neural network structure to avoid laborious calibration of the parameters and to easily incorporate constraints on the solution.

The use of neural networks for solving inverse problems has become increasingly popular, 
especially in the image processing community.
A rich panel of approaches have been proposed, 
either adapted to the sparsity of the  data~\cite{antholzer2019sparse,kofler2018u}, 
or mimicking variational models~\cite{hammernik2018var,adler2017solving}, or iterating learned operators~\cite{aggarwal2018modl,meinhardt2017,pesquet2020learning,hasannasab2020,galinier2020},
or adpating Tikhonov method~\cite{li2020nett}.
% and/or using convolution neural network (CNN) whose parameters are trained on data related to the problem at hand~\cite{antholzer2019sparse,jin2017deep}.
 The successful numerical results of the aforementioned works raise two theoretical questions: when these methods are based on the iteration of a neural network, do they converge (in the sense of the algorithm)? Are these inversion methods stable or robust? %Can we topologically characterize the underlying or implicit regularization of such a method?
 
In iterative approaches,
a regularization operator is learned, either in the form of a proximity (or denoiser) operator as~\cite{meinhardt2017,aggarwal2018modl,galinier2020}, of a regularization term~\cite{li2020nett},
of a pseudodiffential operator~\cite{bubba2021deep}, 
or of its gradient~\cite{gilton2019,Wu2020}. 
Strong connections also exist with Plug and Play methods~\cite{ryu2019plug,sun2019online,pesquet2020learning}, where the regularization operator is a pre-trained neural network.
Such objects have in particular enable high-quality imaging restoration or tomography inversion~\cite{Wu2020}. Here, the non-expensiveness of the neural network is a core property to establish convergence of the algorithm~\cite{Wu2020,pesquet2020learning}.
But our proposed neural network is not based on this idea.

Other recent works solve linear inverse problems by unrolling the optimization iterative process in the form of a network architecture as in~\cite{borgerding2016,jin2017deep}.
Here the number of iterations is fixed, instead of iterating until convergence, and the network is trained in an end-to-end fashion.
Since neural network frameworks offer powerful differential programming capabilities, such architecture are also used for learning hyper-parameters in an unrolled optimization algorithm as in~\cite{Corbineau2020,natarajan2020particle}.

All of the above strategies have shown very good numerical results.
However, few studies have been conducted on their theoretical properties, especially their stability.
The study of the robustness of such a structure is often based on a series of numerical tests, as performed in~\cite{genzel2020robust}.
In~\cite{li2020nett}, they provide very large assumption under which the convergence and the regularization property of their network is ensured. But their result is not subject to verification during the numerical implementation.
A fine characterization of the convergence conditions of recurrent neural network and of their stability via the estimation of a Lipschitz constant is done in~\cite{Combettes2019,Combettes2020}. 
In particular, the Lipschitz constant estimated in~\cite{Combettes2020} is more accurate than in basic approaches 
which often rely in computing the product of the norms of the linear weight operators
of each layer as in~\cite{serrurier2020achieving,cisse2017parseval}.
Thanks to the aforementioned works, proofs of convergence and stability have been demonstrated on specific neural networks applied to inverse problems as in~\cite{pesquet2020learning,hasannasab2020,Corbineau2020}. The analysis carried out in this article is in the line of these references.
%It is in the continuity of this type of results that this work is situated.
% In~\cite{li2020nett} and in~\cite{schwab2019deep}, Tikhonov regularization strategies are developed for an iterative neural network 
% and convergences properties of the network are shown, but little is known about the characterization of the limit point and 
% robustness properties.

Our contributions in this paper are the following. 
\begin{enumerate}
    \item We propose an algorithm based on a neural network architecture to solve the inverse problem~\eqref{def:InversePb}, where a constraint is imposed on the sought solution. One of the main advantages is that the structure of the neural network is interpretable and that it contains few parameters which are learnt.
    \item We study theoretically and numerically the stability of the so-built neural network. The sensitivity analysis is performed with respect to the observed data $y^\delta$ which corresponds to a bias term in each of the layers of \eqref{def:modelNN}.
This analysis is more general than the one performed in \cite{Corbineau2020}, in which only the impact of the initialization was considered.
    \item We show how to implement the neural network in the case of Abel operators. Such operators arise in various physical applications. The proposed 
    neural network performs numerically well compared to other classical inversion methods. Neural network techniques have been widely applied to imaging inverse problems, but few are tested on experimental 1D signal inverse problems.
\end{enumerate}

\paragraph{Outline} 
The outline of the paper is as follows.
In Section~\ref{section:notation}, 
we recall the theoretical background of our work.
We specify our notation, which is based on~\cite{bauschke2011}.
In Section~\ref{section:theory}, 
we establish the stability of the neural network 
defined by~\eqref{def:modelNN},
based on the results in~\cite{Combettes2019} and~\cite{Combettes2020}.
By stability, we mean that the output of the neural network 
is controlled not only with respect to  
its initial input $x_{0}$, but also also with respect to the bias term $T^* y^\delta$. 
The objective is to guarantee that a small difference or error on these vectors is not amplified through the network.
Our theoretical study concerns a class of dynamical systems including a leakage factor,
which is more general  than the neural network defined by~\eqref{def:modelNN}.
In Section~\ref{section:numerics}, the numerical resolution of the problem~\eqref{def:InversePb} is described. 
We define the Abel operator, its characteristics as well as its discretization.
Then, we detail the construction of the training data set.
The architecture of the neural network is explained, 
as well as the sub-network used to estimate the parameters.
Finally, we compute the Lipschitz constants of our trained networks.
We also compare the obtained results with those delivered by two other methods classically used to solve inverse problems involving an Abel integral.

%%%%%%%%%%%%%%%%%%%%%%%%%%%%%%%%%%%%%%%%%%%%%%%%%%%%%%%%%%%%%%%%%%%%%%%%%%%%%%%%%%%%%%%%%%%%%%%%%%%%%%%%%%%%%%%%%%%%%%%%%%%%%%%%%
\section{Notation}
\label{section:notation}

We introduce the theory of convex analysis we will be dealing with, namely monotone operator in Hilbert spaces. We also cover the bits of operator theory that will be needed throughout. 

Let us consider the Hilbert space $\X$ 
endowed with the norm $\| \cdot\|$ and the scalar product $\langle \cdot, \cdot \rangle$. 
In the following, 
$\X$ shall always refer to spaces of functions defined on the interval $]0,1[$.
The notation $\|\cdot\|$ will also refer to the operator norm of bounded operators from $\X$ onto $\X$.
The identity operator over $\X$ will be referred to as $\1$.

An operator $S \colon \X \to \X $ is nonexpansive if it 
is $1-$Lipschitz, that is 
\[
(\forall (x,y) \, \in \, \X \times\X)
\qquad
\| Sx -Sy\| \leq \| x-y\| 
\; .
\]
Moreover, $S$ is said to be
%%%%%%%%%%%%%%%%%
\begin{enumerate}
    \item firmly nonexpansive if
\[
(\forall (x,y) \, \in \, \X \times\X)
\qquad
\| Sx -Sy\|^2
+
\| (\1-S)x -(\1-S)y\|^2
\leq \| x-y\|^2 
\; ;
    \]
    \item a Banach contraction if there exists $\kappa\in ]0,1[$ such that
\begin{equation}\label{e:strictcont}
(\forall (x,y) \, \in \, \X \times\X)
\qquad
\| Sx -Sy\| \le \kappa \| x-y\| 
\; .
\end{equation}
\end{enumerate}
%%%%%%%%%%%%%%%%
If $S$ is a Banach contraction, 
then the iterates $(S^n x)_{n \in \N}$ converge linearly to a fixed point of $S$
according to Picard's theorem.
On the other hand, when $S$ is nonexpansive, 
the convergence is no longer guaranteed.
A way of recovering the convergence of the iterates 
is to assume that  $S$ is averaged,
i.e., 
there exists $\alpha \in ]0,1[$ 
and a nonexpansive operator $R:\X \to \X$ 
such that $S= (1-\alpha) \1+\alpha R$.
In particular, $S$ is $\alpha-$averaged
if and only if
\[
(\forall (x,y) \, \in \, \X \times\X)
\qquad
\| Sx -Sy\|^2
+
\frac{1-\alpha}{\alpha} \| (\1-S)x -(\1-S)y\|^2
\leq \| x-y\|^2 
\; .
\]
If $S$ has a fixed point and it is averaged, 
then the iterates $(S^n x)_{n \in \N}$ converge weakly to a fixed point.
Note that
$S$ is firmly nonexpansive if and only
if it is $1/2-$averaged and that, if $S$ satisfies \eqref{e:strictcont} with $\kappa \in ]0,1[$,
then it is $(\kappa+1)/2$-averaged.

Let $\Gamma_0(\X)$ be the set of proper lower semicontinuous convex function from $\X$ to $]-\infty,+\infty]$.
Then we define the proximal operator as
\begin{definition}
Let $f \in \Gamma_0(\X)$,
$x \in \X$,
and $\gamma>0$.
Then
$\prox_{\gamma f} (x)$ is the unique point that satisfies
\[
\prox_{\gamma f} (x) \,
= \, 
\underset{y \in \X }{\text{argmin}} \, \left(
f(y) + \frac{1}{2\gamma} \| x-y \|^2
\right)
\; .
\]
The function $\prox_{\gamma f} : \X \to \X $ is the proximity operator of $\gamma f$.
\end{definition}
Finally, the proximity operator has the following property.
\begin{proposition}[Proposition 12.28 of~\cite{bauschke2011}]
The operators
$\prox_{\gamma f}$ and $\1-\prox_{\gamma f}$
are firmly nonexpansive.
\end{proposition}

In the proposed neural network~\eqref{def:modelNN}, 
the activation operator is a proximity operator. 
In practice, this is the case for most activation operators, as shown in~\cite{Combettes2019}.
The neural network~\eqref{def:modelNN} is thus a cascade of firmly nonexpansive operators and linear operators. 
If the linear part is also nonexpansive, bounds on the effect of a pertubation of the neural network or its iterates can be established.

%%%%%%%%%%%%%%%%%%%%%%%%%%%%%%%%%%%%%%%%%%%%%%%%%%%%%%%%%%%%%%%%%%%%%%%%%%%%%%%%%%%%%%%%%%%%%%%%%%%%%%%%%%%%%%%%%%%%%%%%%%%%%%%%%
\section{Stability and $\alpha$-averagedness}
\label{section:theory}
%
%%%%%%%%%%%%%%%%%%%%%%%%%%%%%%%%%%%%%%%%%%
In this section, we study the stability of the proposed neural network~\eqref{def:modelNN}. 
This analysis is performed by estimating the Lipschitz constant of the network, and by determining under which conditions this network is $\alpha$-averaged. 
To do so, we introduce a virtual network, which takes as inputs the classical ones on top of a new one, which is the bias parameter.

\subsection{Virtual neural network with leakage factor}\label{se:VNN1}
To facilitate our theoretical analysis, we will introduce a virtual network making use of new variables $(z_n)_{n\in \mathbb{N}}$. 
For every $n \in \mathbb{N}\setminus\{0\}$,
we define the $n$-th layer of our virtual network as follows
\begin{equation}
    \label{def:nn-virtual-bias}
%    (\forall n \in \mathbb{N})\quad
z_n = \left( \begin{array}{c}
 x_n \\
 b_n
\end{array}
\right) ,
\quad
z_{n} = Q_n(U_n z_{n-1})
\; ,
\quad
\text{with}
\quad
\begin{cases}
\displaystyle Q_n = \left( \begin{array}{c}
  R_n   \\
  \1
\end{array} \right)
\; , \\
\\
\displaystyle U_n = \left( \begin{array}{cc}
   W_n  &  \lambda_n \1 \\
   0    & \eta_n \1
\end{array}\right)
\;.
\end{cases}
\end{equation}
%where $\1$ denotes the identity operator.
Note that, in order to gain more flexibility, we have 
 included positive multiplicative factors $(\eta_n)_{n\geq 1}$ on the bias.
 Cascading $m$ such layers yields
  \begin{equation}
	\label{def:modelNN-leakage}
	\begin{cases}
	\textbf{Initialization:} \\
	\quad b_0 = T^*\yd ,\\
	\textbf{Layer $n\in \{1,\ldots,m\}$:} \\
      \quad x_n %=  Q_n(x_0,x_{n-1}) 
          = R_n(W_n x_{n-1} + V_n b_0)\;,
    \end{cases}		  
\end{equation}
where 
\begin{align}
&R_n = \text{prox}_{\lambda_n \mu_n g}\\
&W_n =  \1 - \lambda_n T^*T - \lambda_n \tau_n D^*D \\
&V_n = \lambda_n \eta_{n-1}\cdots \eta_{1} \1
\end{align}
and $\eta_{0}=1$.
We thus see that the  network defined by Model~\eqref{def:modelNN} 
is equivalent to the virtual one when all the factors $\eta_{n}$ are equal to one.
When $n\ge 1$ and $\eta_{n}< 1$.
The  parameters $(\eta_{n})_{n\geq 1}$
can be interpreted as a leakage factor. 
\begin{remarque}
In the original forward-backward algorithm, the introduction of $(\eta_{n})_{n\geq 1}$ amounts to introducing an error $e_n$ in the gradient step, at iteration $n$, which is equal to
\begin{equation}
    e_n = \lambda_n (\eta_{n-1}\cdots \eta_1-1) b_0.
\end{equation}
From known properties 
concerning the forward-backward algorithm~\cite{combettes2011proximal}, 
the convergence of the algorithm is still guaranteed provided that
\begin{equation}
    \sum_{n=2}^{+\infty} \lambda_n |\eta_{n-1}\cdots \eta_1-1| < +\infty
    \; .
\end{equation}
\end{remarque}
In our analysis, it will be useful to define the triangular linear operator
\begin{equation}
    \label{def:U-bias}
U = \;U_m \circ \cdots \circ U_1 
  = \;
\left( \begin{array}{cc}
  W_{1,m}   &  \widetilde{W}_{1,m}\\
  0   &  \eta_{1,m} \1
\end{array} \right)
\; ,
\end{equation}
where, for every  $n\in \{1,\ldots,m\}$ and $i\in \{1,\ldots,n\}$ 
\begin{equation}\label{e:defWtin}
\widetilde{W}_{i,n} =  
\sum_{j=i}^{n} \lambda_{j} \eta_{i,j-1}W_{j+1,n}\, 
\end{equation}
%\begin{align}
%\widetilde{W} &=   
%  \sum_{n=1}^{m}  \lambda_{n} \eta_{1,n-1} W_{n+1,m} 
%  \end{align}
and, for every $i\in \{1,\ldots,m+1\}$
and  $j\in \{0,\ldots,m\}$,
\begin{align}
&   W_{i,j} = 
\begin{cases}
W_j\circ \cdots \circ W_i & \mbox{if $j\ge i$}\\
\1 & \mbox{otherwise,}
\end{cases}\label{def:WjWi}\\
&    \eta_{i,j} = 
    \begin{cases}
    \eta_{j}\cdots \eta_{i} & \mbox{if $j\ge i$}\\
    1 & \mbox{otherwise.}
    \end{cases}\label{def:N-bias}
\end{align}
Since $T$ defined by~\eqref{def:T} is a compact operator, 
we can define its singular value expansion as in~\cite{Engl1996}.
Furthermore, we place ourselves in the case where
$D^*D$ and $T^*T$ commutes, 
for operators $T$ defined by~\eqref{def:T} and regularization operators $D$. 
Therefore those operators admit the same eigensystem.
In particular, they can be diagonalized in the same orthonormal set of eigenvectors $(v_{p})_{p}$.
We define their respective eigenvalues
$(\beta_{T,p})_p$ and $(\beta_{D,p})_p$, as well as 
the following quantities, 
for every eigenspaces $p\in \mathbb{N}$, 
$n\in \{1,\ldots,m\}$, and $i\in \{1,\ldots,n\}$,
\begin{flalign}
\label{def:vp0}
     &\beta_p^{(n)} = 1 - \lambda_n (\beta_{T,p} +\tau_n \beta_{D,p}) \; , &\\
\label{def:vp1}
     &\beta_{i,n,p} = \prod_{j=i}^n \beta_p^{(j)} \; , & \\ 
\label{def:vp-bias}
     &  \widetilde{\beta}_{i,n,p} =   \sum_{j=i}^{n-1} \beta_p^{(n)} \cdots \beta_p^{(j+1)} \lambda_{j}\eta_{i,j-1}
     + \lambda_n \eta_{n-1}\cdots \eta_i
     \;  & 
\end{flalign}
with the convention $\sum_{i=n}^{n-1} \cdot = 0$.
Note that $(\beta_{i,n,p},v_p)_{p\in\mathbb{N}}$
and $(\widetilde{\beta}_{i,n,p},v_p)_{p\in\mathbb{N}}$
are the eigensystems of~$W_{i,n}$ and~$\widetilde{W}_{i,n}$, respectively.

\subsection{Stability results for the virtual network}
We first recall some recent results on the stability of neural networks \cite[Proposition 3.6(iii)]{Combettes2019} \cite[Theorem 4.2]{Combettes2020}.
\begin{proposition}\label{prop:oldres}
Let $m > 1$ be an integer, let $(\mathcal{H}_i)$ be nonzero real Hilbert spaces.
For every $n \in \{1,\ldots,m\}$, let $U_n \in \mathcal{B}(\mathcal{H}_{n-1},\mathcal{H}_{n})$
and let $Q_n\colon \mathcal{H}_n \to \mathcal{H}_n$ be a firmly nonexpansive operator.
Set $U = U_{m} \circ \cdots \circ U_1$ and
\begin{multline}
\label{e:defthetaell}
\theta_m=\|U\|\\
+\sum_{k=1}^{m-1}\sum_{1\leq j_1<\ldots<j_k\leq m-1}
\|U_m\circ\cdots\circ U_{j_k+1}\|\,
\|U_{j_k}\circ\cdots\circ U_{j_{k-1}+1}\|\cdots 
\|U_{j_1}\circ\cdots\circ U_1\|.
\end{multline}
Let $S = Q_m\circ U_m \circ \cdots \circ Q_1 \circ U_1$.
Then the following hold
\begin{enumerate}
    \item\label{prop:oldresi} $\theta_m/2^{m-1}$ is a Lipschitz constant of $S$.
    \item\label{prop:oldresii} 
    Let $\alpha \in [1/2, 1]$.
    If $\mathcal{H}_m = \mathcal{H}_0$ and 
    \begin{equation}\label{e:condalphaold}
\| U - 2^{m} (1- \alpha) \1 \|
- \|U\|
+2 \theta_{m}
\leq
2^{m} \alpha
\;,
\end{equation}
then $S$ is $\alpha$-averaged.
\end{enumerate}
\end{proposition}
In light of these results, 
we will now analyze the properties of the virtual network~\eqref{def:nn-virtual-bias}
based on the singular values of the operators
$T$ and $D$, 
%defined by~\eqref{def:T} and~\eqref{def:D}, 
and the parameters $(\lambda_n)_{1\le n\le m}$ and $(\tau_n)_{1\le n \le m}$. 
One of the main difficulties with respect to the case already studied by~\cite{Corbineau2020}
is that here the involved operators $(U_n)_{1\le n \le m}$ are no longer self-adjoint.

A preliminary result will be needed:
\begin{lemme}\label{e:normUni-bias}
Let $m \in \mathbb{N}\setminus \{0\}$
the total number of layers.
For every layer $n\in \{1,\ldots,m\}$ and 
layer $i\in \{1,\ldots,n\}$, 
the norm of $U_n\circ \cdots \circ U_i$ is 
equal to $\sqrt{a_{i,n}}$ with
\begin{equation}\label{e:defain-bias}
a_{i,n} = \frac{1}{2} \sup_{p \in \N}
\left(
\beta_{i,n,p}^2 + \widetilde{\beta}_{i,n,p}^2 + \eta_{i,n}^2+ \sqrt{(\beta_{i,n,p}^2  + \widetilde{\beta}_{i,n,p}^2  +  \eta_{i,n}^2 )^2
        - 4\beta_{i,n,p}^2 \eta_{i,n}^2}
\right) \; ,
\end{equation}
where $p$ covers the eigenspaces of $T^*T$ defined by~\eqref{def:T}.
\end{lemme}
\begin{preuve}
Thanks to expressions \eqref{def:nn-virtual-bias}, \eqref{e:defWtin}, \eqref{def:WjWi}, and \eqref{def:N-bias},
%~\eqref{def:vp1}, \eqref{def:N-bias} and~\eqref{def:vp-bias}, 
we can calculate the norm of $\| U_n\circ \cdots \circ U_i \|$.
For every $z =(x,b)$, $U_n\circ \cdots \circ U_i z = (W_{i,n} x + \widetilde{W}_{i,n} b, \eta_{i,n} b)$ and
\[
\begin{split}
\|U_n\circ \cdots \circ U_i z \|^2 = & \; \| W_{i,n} x + \widetilde{W}_{i,n} b\|^2 + \eta_{i,n}^2 \|b\|^2\\
= & \; \| W_{i,n}x\|^2 + 2\langle W_{i,n} x , \widetilde{W}_{i,n} b \rangle 
+ \| \widetilde{W}_{i,n} b \|^2 + \eta_{i,n}^2 \|b\|^2
\; .
\end{split}
\]
%where 
%\begin{equation}
%\widetilde{W}_{i,n} =  
%\sum_{j=i}^{n} \lambda_{j} \eta_{i,j-1}W_{j+1,n}\, .
%\end{equation}
Let $(\beta_{i,n,p},v_p)_{p\in\mathbb{N}}$ 
defined by~\eqref{def:vp1}
and $(\widetilde{\beta}_{i,n,p},v_p)_{p\in\mathbb{N}}$
defined by~\eqref{def:vp-bias}
be the respective eigensystems of~$W_{i,n}$ and~$\widetilde{W}_{i,n}$.
Let us decompose $(x,b)$ in the basis of eigenvectors $(v_{p})_{p}$ of $T^*T$, % defined by~\eqref{def:T},
as 
\[
\begin{cases}
x = \sum_p \xp \, v_p \; , \\
b = \sum_p \bp \, v_p \; .
\end{cases} 
\]
We have then
\[
\|U_n\circ \cdots \circ U_i z\|^2 = \; \sum_p \beta_{i,n,p}^2 \xp^2 + 2 \sum_p \beta_{i,n,p} \widetilde{\beta}_{i,n,p} \xp \bp
+ \sum_p (\widetilde{\beta}_{i,n,p}^2 + \eta_{i,n}^2) \bp^2
\; .
\]
By definition of the operator norm,
\[
 \| U_n\circ \cdots \circ U_i\|^2 = \; \underset{\|x\|^2 +\|b\|^2=1}{\sup}
\left(
\sum_p \beta_{i,n,p}^2 \, \xp^2 
+ (\eta_{i,n}^2+\widetilde{\beta}_{i,n,p}^2) \, \bp^2 
+ 2\beta_{i,n,p} \widetilde{\beta}_{i,n,p} \, \xp  \,\bp \right)\; .
\]
Note that, for every integer $p \in \N$ and $\zp = (\xp, \bp)\in \mathbb{R}^2$,
\begin{equation}
\beta_{i,n,p}^2 \, \xp^2 +
(\eta_{i,n}^2+\widetilde{\beta}_{i,n,p}^2) \, \bp^2 
+ 2\beta_{i,n,p} \widetilde{\beta}_{i,n,p} \, \xp \, \bp 
= 
\langle A_{i,n,p} \zp, \zp \rangle
\end{equation}
where $\langle \cdot, \cdot \rangle$ denotes the Euclidean inner product and 
\[A_{i,n,p} = \left( 
\begin{array}{cc}
  \beta_{i,n,p}^2   & \beta_{i,n,p} \widetilde{\beta}_{m,p} \\
  \widetilde{\beta}_{i,n,p}\beta_{i,n,p}     & \eta_{i,n}^2 + \widetilde{\beta}_{i,n,p}^2
\end{array}
\right) \; .
\]
Hence,
\[
 \| U_n\circ \cdots \circ U_i \|^2 = \; \underset{z = (\zp)_p, \|z\| =1}{\sup}
             \sum_p \langle A_{i,n,p} \, \zp, \, \zp \rangle 
\; .
\]
Since $A_{i,n,p}$ is a symmetric positive semidefinite matrix, 
\begin{equation}\label{e:defnormU2-bias}
\| U_n\circ \cdots \circ U_i \|^2 = \underset{p, \|\zp\| =1}{\sup} \langle A_{i,n,p} \zp, \zp \rangle = \sup_p\, \nu_{i,n,p}\; ,
\end{equation}
where, for every $p\in \mathbb{N}$, $\nu_{i,n,p}$ is  the maximum eigenvalue of $A_{i,n,p}$.
The two eigenvalues of this matrix are the roots of the characteristic polynomial 
\[
\begin{split}
(\forall \nu \in \mathbb{R})\quad
    \text{det} (A_{i,n,p} - \nu \1_2) = & \;
(\beta_{i,n,p}^2 - \nu)(\widetilde{\beta}_{i,n,p}^2 + \eta_{i,n}^2 - \nu) - \beta_{i,n,p}^2\widetilde{\beta}_{i,n,p}^2 \\
                           = & \;
\nu^2 - (\beta_{i,n,p}^2 + \widetilde{\beta}_{i,n,p}^2 + \eta_{i,n}^2 )\nu 
+ \beta_{i,n,p}^2\eta_{i,n}^2
\; .
\end{split}
\]
The discriminant of this second-order polynomial reads
\[
\begin{split}
\Delta_{i,n,p} = & \;  (\beta_{i,n,p}^2 + \widetilde{\beta}_{i,n,p}^2 + \eta_{i,n}^2 )^2 - 4\beta_{i,n,p}^2 \eta_{i,n}^2\\
        =& \;  (\beta_{i,n,p}^2 -\widetilde{\beta}_{i,n,p}^2 - \eta_{i,n}^2)^2
       +4 \beta_{i,n,p}^2\widetilde{\beta}_{i,n,p}^2 \; \ge 0
        \; .
\end{split}
\]
Therefore, for every $p\in \mathbb{N}$,
\begin{equation}\label{e:defnuinp-bias}
\nu_{i,n,p} = \frac{1}{2} 
\left(
\beta_{i,n,p}^2 + \widetilde{\beta}_{i,n,p}^2 + \eta_{i,n}^2 + \sqrt{(\beta_{i,n,p}^2  + \widetilde{\beta}_{i,n,p}^2  + \eta_{i,n}^2 )^2
        - 4\beta_{i,n,p}^2 \eta_{i,n}^2}
\right) \; .
\end{equation}
By going back to \eqref{e:defnormU2-bias},
we obtain
\[
\| U_n\circ \cdots \circ U_i  \|^2 = a_{i,n} \; .
\]
\end{preuve}

We will now quantify the Lipschitz regularity of the network.
\begin{proposition}\label{p:LipVNN-bias}
Let $m \in \mathbb{N}\setminus \{0\}$. 
For every $n\in \{1,\ldots,m\}$ and 
$i\in \{1,\ldots,n\}$,  let $a_{i,n}$ be given by
\eqref{e:defain-bias}.
Set $\theta_0 = 1$ and
define $(\theta_n)_{1\le n \le m}$ recursively by
\[
(\forall n \in \{1,\ldots,m\})\quad 
\theta_n = 
\sum_{i=1}^n \theta_{i-1} \sqrt{a_{i,n}}\;.
\]
Then $\theta_m/2^{m-1}$ is a Lipschitz constant of the virtual network~\eqref{def:nn-virtual-bias}.
\end{proposition}
\begin{preuve}
According to Proposition~\ref{prop:oldres}\ref{prop:oldresi}, 
if $\theta_m$ is given by \eqref{e:defthetaell}, then $\theta_m/2^{m-1}$ is a Lipschitz constant of the virtual network~\eqref{def:nn-virtual-bias}. 
On the other hand, it follows from
\cite[Lemma 3.3]{Combettes2019} that $\theta_m$ can be calculated recursively as
\[
(\forall n \in \{1,\ldots,m\})\quad 
\theta_n = 
\sum_{i=1}^n \theta_{i-1}  \| U_n \circ \cdots \circ U_i \|\;,
\]
with  $\theta_0 = 1$. Finally, Lemma \eqref{e:normUni-bias} allows us to substitute $(\sqrt{a_{i,n}})_{1\le i \le n}$ for $(\| U_n \circ \cdots \circ U_i \|)_{1\le i \le n}$ in the above expression.
\end{preuve}

We will next provide conditions ensuring that the virtual network is an averaged operator.
\begin{proposition}
\label{prop:ab-bias}
Let $m \in \mathbb{N}\setminus \{0\}$.
Let $a_{1,m}$ be defined in Lemma \ref{e:normUni-bias} and $\theta_m$ be defined in Proposition \ref{p:LipVNN-bias}.
Let $\alpha \in [1/2,1]$.
Define
\begin{align}
b_\alpha =  \frac{1}{2} \sup_p &\;
\left(
 (\beta_{1,m,p}- \gamma_\alpha)^2  +  (\eta_{1,m}- \gamma_\alpha)^2 +\widetilde{\beta}_{1,m,p}^2 \right. \nonumber\\
& \; \left. + \; \sqrt{
\begin{array}{c}
\big( (\beta_{1,m,p}- \gamma_\alpha)^2  +  (\eta_{1,m}-\gamma_\alpha)^2 + \widetilde{\beta}_{1,m,p}^2 \big)^2 \\
 - 4  (\beta_{1,m,p}- \gamma_\alpha)^2(\eta_{1,m}- \gamma_\alpha)^2
\end{array}
}
\; \right) \; ,
\end{align} 
with $\gamma_\alpha = 2^m(1-\alpha)$.
Then virtual network~\eqref{def:nn-virtual-bias} is $\alpha$-averaged if
\begin{equation}\label{e:condalphanew-bias}
\sqrt{b_\alpha} -\sqrt{a_{1,m}} \leq 2^m \alpha - 2\theta_m\;.
\end{equation}
\end{proposition}

\begin{preuve}
Let us calculate
the operator norms of $ U $ and 
$ U - \gamma_{\alpha} \1 $,
where $U$ is given by \eqref{def:U-bias}.

\paragraph{Norm of $U$. }
Applying Lemma~\ref{e:normUni-bias} when $i=1$ and $n=m$ yields
\[
\| U\|^2 = a_{1,m} 
\; .
\]

\paragraph{Norm of $U - \gamma_{\alpha} \1$. }
We follow the same reasoning as in the proof of Lemma~\ref{e:normUni-bias}.
We have
\[
\begin{split}
\| U -\gamma_{\alpha} \1 \|^2 
                        =  \underset{z=(\zp)_p, \|z\| =1}{\sup} 
\sum_p \langle B_p \zp , \zp \rangle
\; ,
\end{split}
\]
where $B_p$ is the symmetric positive semidefinite matrix given by
\[B_p = \left(
\begin{array}{cc}
 (\beta_{1,m,p}-\gamma_{\alpha})^2    & (\beta_{1,m,p} -\gamma_{\alpha})\widetilde{\beta}_{1,m,p}  \\
  (\beta_{1,m,p} -\gamma_{\alpha})\widetilde{\beta}_{1,m,p}     &  (\eta_{1,m}-\gamma_{\alpha})^2 +\widetilde{\beta}_{1,m,p}^2
\end{array}
\right) \; .
\]
By definition of the spectral norm, 
\begin{equation}\label{e:norm2Ualpha-bias}
\| U -2^m(1-\alpha) \1 \|^2  = \sup_p\, \nu_{p}\; ,
\end{equation}
where, for every $p\in \mathbb{N}$, $\nu_{p}$ is  the maximum eigenvalue of $B_{p}$.
The two eigenvalues of this matrix are the roots of the polynomial 
\begin{align}
(\forall \nu \in \mathbb{R})\quad
    \operatorname{det} (B_{p} - \nu \1_2) 
                           = & \;
 \nu^2 - ( (\beta_{1,m,p}-\gamma_\alpha)^2  +  (\eta_{1,m}-\gamma_\alpha)^2 +\widetilde{\beta}_{1,m,p}^2 ) \ \nu 
\nonumber\\&+(\beta_{1,m,p}-\gamma_\alpha)^2 (\eta_{1,m}-\gamma_\alpha)^2
\; .
\end{align}
Solving the corresponding second-order equation leads to
\begin{align}
\sup_p \; \nu_p =  &\frac{1}{2} 
\bigg(
 (\beta_{1,m,p}-\gamma_\alpha)^2  +  (\eta_{1,m}-\gamma_\alpha)^2 +\widetilde{\beta}_{1,m,p}^2 \nonumber\\
& \; + \sqrt{\big( (\beta_{1,m,p}-\gamma_\alpha)^2  +  (\eta_{1,m}-\gamma_\alpha)^2 + \widetilde{\beta}_{1,m,p}^2 \big)^2 - 4  (\beta_{1,m,p}-\gamma_\alpha)^2(\eta_{1,m}-\gamma_\alpha)^2}
\bigg) \; .
\end{align}
Then, it follows from~\eqref{e:norm2Ualpha-bias} that $\|U-\gamma_{\alpha}\1\|^2 = b_\alpha$.

\paragraph{Conclusion of the proof.} 
Based on the previous calculations, Condition \eqref{e:condalphanew-bias} is equivalent to \eqref{e:condalphaold}.
In addition, let us note that for every $n\in \{1,\ldots,m\}$, $Q_n$ 
in \eqref{def:nn-virtual-bias} is firmly nonexpansive since $R_n$ and $\1$ are.
By applying now Proposition~\ref{prop:oldres}\ref{prop:oldresii},
we deduce that, when Condition \eqref{e:condalphanew-bias} holds, 
 virtual network~\eqref{def:nn-virtual-bias} is $\alpha$-averaged.
\end{preuve}
\begin{remarque}
Condition \eqref{e:condalphanew-bias} just provides a sufficient condition for the averagedness of virtual network~\eqref{def:nn-virtual-bias}.
\end{remarque}
%
%%%%%%%%%%%%%%%%%%%%%%%%%%%%%%%%%%%%%%%%%%
\subsection{Link with the original neural network -- direct approach}
In this subsection we go back to our initial model defined by~\eqref{def:modelNN-leakage}.
We consider two different inputs $z_1=(x_1,b_1)$ and $z_2=(x_2,b_2)$ in $\X \times \X$. The distance between these points is
\[
\|z_2-z_1\| = \sqrt{\|x_2-x_1\|^2 + \|b_2-b_1\|^2}
\; .
\]
Let $z_{i,n}=(x_{i,m},b_{i,m})$ be the output of the $m$-th layer of virtual network~\eqref{def:nn-virtual-bias}.
Then, 
\[
\begin{split}
 \| z_{2,m} - z_{1,m} \|^2 = & \; \| x_{2,m} - x_{1,m}\|^2 + \| b_{2,m} - b_{1,m}\|^2
\\
=& \; \| x_{2,m} - x_{1,m}\|^2 + \eta_{1,m}^2 \| b_{2} - b_{1}\|^2 
\; ,\\
\end{split}
\]
and, thanks to Proposition~\ref{p:LipVNN-bias},
\[
\begin{split}
 \| z_{2,m} - z_{1,m} \|^2 \leq \;
 \frac{\theta_m^2}{2^{2(m-1)}}\left(
 \| x_2-x_1\|^2 + \| b_2-b_1 \|^2
 \right)
 \; .\\
\end{split}
\]
Then, the following inequality allows us to quantify the Lipschitz properties of the neural 
network~\eqref{def:modelNN} with an error on $b_n$:
\[
\| x_{1,m} - x_{2,m} \|^2 \leq \;
 \frac{\theta_m^2}{2^{2(m-1)}}
 \| x_2-x_1\|^2 
 + \left( \frac{\theta_m^2}{2^{2(m-1)}} - \eta_{1,m}^2 \right)\| b_2-b_1 \|^2
 \; .
\]
Two cases are of interest:
\begin{itemize}
    \item If the network is initialized with a fixed signal, say $x_{1,0} = x_{2,0}= 0$, then
 \[
 \| x_{1,m} - x_{2,m} \|^2 \leq \; 
\left( \frac{\theta_m^2}{2^{2(m-1)}} - \eta_{1,m}^2 \right) 
\; \| b_2-b_1\|^2
\; .
\]
So, a Lipschitz constant with respect to the 
input data $T^* y^\delta$ is 
\begin{equation}\label{e:Lipreal1-bias}
\vartheta_m = \sqrt{\frac{\theta_m^2}{2^{2(m-1)}} - \eta_{1,m}^2}.
\end{equation}
\item On the other hand, if the initialization is dependent
on the observed image, i.e. $x_{1,0}= b_1$ and $x_{2,0} = b_2$,
\[
\| x_{1,m} - x_{2,m} \|^2 \leq \; 
\left(\frac{\theta_m^2}{2^{2m-3}} -\eta_{1,m}^2\right) 
\; \| b_2-b_1\|^2
\; .
\]
So a higher Lipschitz constant value w.r.t. to the input data is obtained:
    \begin{equation}
\vartheta_m = \sqrt{\frac{\theta_m^2}{2^{2m-3}} -\eta_{1,m}^2}.
\end{equation}
\end{itemize}
\begin{remarque}\label{re:directpasbon}
Let us go back to Model~\eqref{def:modelNN}.
We thus consider the virtual Model~\eqref{def:nn-virtual-bias} 
without leakage factor, i.e.,  for every $n \in \{ 1, \ldots, m \}$, $\eta_n =1$ and $\eta_{1,m}=1$.
Then, for every $n \in \{ 1, \ldots, m \}$,
\begin{equation}
    \label{def:nn-virtual}
z_n = \left( \begin{array}{c}
 x_n \\
 b_n
\end{array}
\right) ,
\quad
z_{n} = Q_n(U_n z_{n-1})
\; ,
\quad
\text{with}
\quad
\begin{cases}
\displaystyle Q_n = \left( \begin{array}{c}
  R_n   \\
  \1
\end{array} \right)
\; , \\
\\
\displaystyle U_n = \left( \begin{array}{cc}
   W_n  &  V_n \\
   0    & \1
\end{array}\right)
\; .
\end{cases}
\end{equation}
 Assume that the virtual network in~\eqref{def:nn-virtual} is $\alpha$-averaged.
Then it is 1-Lipschitz. This is also consistent with
\eqref{e:condalphaold} which implies that 
\begin{equation}
\|U - 2^{m} (1- \alpha)\1 \| - \|U\|+2 \theta_{m} 
\leq
2^{m} \alpha
\quad \Rightarrow \quad \frac{\theta_m}{2^{m-1}} \le 1
\;.
\end{equation}
Then, according to \eqref{e:Lipreal1-bias}, we would get
$\vartheta_m=0$, which would mean that the network delivers an output independent of the available data.
If we except trivial cases for which $R_{n}= 0$ or $W_{n}= 0$ for some $n\in \{1,\ldots,m\}$, this behavior is impossible.
So this means that virtual network~\eqref{def:nn-virtual} cannot be $\alpha$-averaged, 
hence Condition~\eqref{e:condalphanew-bias} is not met when,  for every $n \in \{ 1, \ldots, m \}$,  $\eta_n=1$.

This can be concluded more directly.
For every $\alpha \in ]0,1[$,
virtual network~\eqref{def:nn-virtual} cannot be $\alpha$-averaged, 
since
$R = (1-1/\alpha) \1 + 1/\alpha S$ cannot be nonexpansive.
Indeed suppose that $ \|R(x_1,b_1) - R(x_2,b_2)\|^2 \leq \|x_1 -x_2\|^2 + \|b_1 -b_2\|^2 $, for every $(x_1,b_1)$ and $(x_2,b_2)$ in $\X \times \X$. Since
\[
\begin{split}
\|R(x_1,b_1) - R(x_2,b_2)\|^2 
  = & \; \Big\| (1-1/\alpha) \left( \begin{array}{c}
      x_1-x_2  \\
      b_1-b_2
\end{array} \right)
 + 1/\alpha \left( \begin{array}{c}
      x_{1,m}-x_{2,m}  \\
      b_1-b_2
\end{array} \right) \Big\|^2 \\
  = & \; \|(1-1/\alpha) (x_1-x_2) + 1/\alpha (x_{1,m}-x_{2,m}) \|^2 +\| b_1 -b_2 \|^2 \\
  \leq & \; \|x_1 -x_2\|^2 + \|b_1 -b_2\|^2
  \; ,
\end{split}
\]
we deduce that
\[
 \|(1-1/\alpha) (x_1-x_2) + 1/\alpha (x_{1,m}-x_{2,m}) \|
 \leq 
 \|x_1 -x_2\|
 \; ,
\]
which cannot stand since, for $b_1\neq b_2$, $x_{1,m}-x_{2,m}$ can be nonzero when $x_1=x_2$. 
\end{remarque}

\subsection{Link with the original neural network -- use of a semi-norm}
Remark \ref{re:directpasbon} suggests that we need a finer strategy to evaluate 
the nonexpansiveness properties of
Model~\eqref{def:modelNN}.
%We consider Virtual Model~\eqref{def:nn-virtual}, where $\eta_n =1$
%for $n \in \{ 1, \ldots, m \} $.
On the product space $\X \times \X$, we define the semi-norm which takes only into account the first component of the vectors:
\begin{equation}
\label{def:seminormv}
    z = (x,b) \mapsto |z| = \|x\|.
\end{equation}
Let $L\colon \X\times \X \to \X\times \X$ be any bounded linear operator and, for every $z\in \X \times \X$, 
let $Lz = ((Lz)_{\rm x},(Lz)_{\rm b})$.
We define the associated operator semi-norm
\begin{equation}
    \label{def:seminorm}
| L | = \underset{\|z\| = 1 }{\text{sup}} \| (Lz)_{\rm x}\|
\; .
\end{equation}
\begin{lemme}\label{e:normUni-proj}
Let $m \in \mathbb{N}\setminus \{0\}$.
For every $n\in \{1,\ldots,m\}$ and 
$i\in \{1,\ldots,n\}$, 
the seminorm  $|U_n\circ \cdots \circ U_i|$ is 
equal to $\sqrt{\overline{a}_{i,n}}$ with
\begin{equation}\label{e:defain-proj}
\overline{a}_{i,n} = \sup_p
\left(
\beta_{i,n,p}^2 + \widetilde{\beta}_{i,n,p}^2
\right) \; .
\end{equation}
\end{lemme}
\begin{preuve}
The seminorm of $U_n\circ \cdots \circ U_i$ is the same as the norm
of $U_n\circ \cdots \circ U_i$ where $\eta_{n}$ has been set to 0.
The result thus follows from Lemma \ref{e:normUni-bias} where $\eta_{i,n}= 0$.
\end{preuve}
\begin{proposition}
\label{p:LipVNN-proj}
Let $m \in \mathbb{N}\setminus \{0\}$. 
For every $i\in \{1,\ldots,n\}$ and $n\in \{1,\ldots,m-1\}$,
let $a_{i,n}$ be defined by
\eqref{e:defain-bias} and
let $\overline{a}_{i,m}$ be given by
\eqref{e:defain-proj}.
Set $\theta_0 = 1$ and
define 
\begin{align}
&(\forall n \in \{1,\ldots,m-1\})
\quad 
\theta_n = 
\sum_{i=1}^n \theta_{i-1} \sqrt{a_{i,n}}\;,\\
&\overline{\theta}_m = 
\sum_{i=1}^m \theta_{i-1} \sqrt{\overline{a}_{i,m}}\;.
\end{align}
Then
the network in~\eqref{def:modelNN-leakage} with input $(x_0,b_0)$ and output $x_m$ is $\overline{\theta}_m /2^{m-1}$-Lipschitz.
\end{proposition}
\begin{preuve}
Network~\eqref{def:modelNN-leakage}
can be expressed as $R_m \circ \overline{U}_m \circ Q_{m-1}\circ U_{m-1}\circ \cdots \circ Q_1 \circ U_1$ where 
\begin{equation}\label{e:oUm}
\overline{U}_m = D_{\rm x}\circ U_m
\end{equation}
and $D_{\rm x}$ is the decimation operator
\begin{equation}
    D_{\rm x}= [\1\quad 0].
\end{equation}
This  network has the same Lipschitz properties as the network in~\eqref{def:nn-virtual-bias} with $\eta_{m}= 0$.
The result can thus be deduced from Proposition \ref{p:LipVNN-bias} by setting $\eta_{1,m}= 0$.
\end{preuve}
To investigate averagedness properties, a first possibility is to consider a network from $\X\times \X$ to $\X\times \X$ with input $(x_{0},b_{0})$ and
output $(x_{m},0)$. 

\begin{proposition}
\label{prop:ab-proj}
Let $m \in \mathbb{N}\setminus \{0,1\}$. 
Let $\overline{a}_{1,m}$ be defined in Lemma~\ref{e:normUni-proj} and $\overline{\theta}_m$ be defined in Proposition \ref{p:LipVNN-proj}.
Let $\alpha \in [1/2,1]$.
Define
\begin{align}
\overline{b}_\alpha =  \frac{1}{2} \sup_p &\;
\left(
 (\beta_{1,m,p}- \gamma_\alpha)^2  +\widetilde{\beta}_{1,m,p}^2+  \gamma_\alpha^2 + \right. \nonumber\\
& \; \left. + \; \sqrt{
\begin{array}{c}
\big( (\beta_{1,m,p}- \gamma_\alpha)^2  + \widetilde{\beta}_{1,m,p}^2+  \gamma_\alpha^2  \big)^2
 - 4  (\beta_{1,m,p}- \gamma_\alpha)^2\gamma_\alpha^2
\end{array}
}
\; \right) \; ,
\end{align} 
with $\gamma_\alpha = 2^m(1-\alpha)$.
 If
 \begin{equation}\label{e:condalphanew-proj}
\sqrt{\overline{b}_\alpha} -\sqrt{\overline{a}_{1,m}} \leq 2^m \alpha - 2\overline{\theta}_m\;,
\end{equation}
then 
the network in~\eqref{def:modelNN} with input $(x_0,b_0)$ and output $(x_m,0)$ is $\alpha$-averaged. 
\end{proposition}

\begin{preuve}
The network of interest is
\[
Q_{m}\begin{bmatrix}
\overline{U}_m \circ Q_{m-1}\circ U_{m-1}\circ Q_1 \circ U_1\\
0
\end{bmatrix}\;.
\]
This network can be viewed as a special case of the network in~\eqref{def:nn-virtual-bias} where $\eta_{m}=0$,
which implies that $\eta_{1,m}= 0$. The result is thus a consequence
of Proposition \ref{prop:ab-bias}.
\end{preuve}

Another possibility for investigating averagedness properties consists of defining a network from $\X$ to $\X$. 
We will focus on two specific networks of the form
\begin{equation}\label{e:NNreal}
R_m\circ \overline{U}_m \circ Q_{m-1}\circ U_{m-1}\cdots
Q_1 \circ \widehat{U}_{1},
\end{equation}
where $\overline{U}_{m}$ is given by \eqref{e:oUm}.
\begin{enumerate}
\item\label{caseNN1} The first one assumes that $x_0 = 0$ in \eqref{def:modelNN}. It is thus given by 
\begin{equation}
\widehat{U}_{1} = U_1 \begin{bmatrix}
0\\
\1
\end{bmatrix}\,.
\end{equation}
\item\label{caseNN2} The second one assumes that $x_0=b_0$ in \eqref{def:modelNN}. It is thus given by
\begin{equation}
\widehat{U}_{1} = U_1 \begin{bmatrix}
\1\\
\1
\end{bmatrix}\,.
\end{equation}
\end{enumerate}

By proceeding similarly to the proof of Lemma \ref{e:normUni-bias}
and Proposition \ref{prop:ab-bias}, we obtain the following two results: 
\begin{lemme}\label{e:normUni-proj2}
Let $m \in \mathbb{N}\setminus \{0,1\}$.
For every $n\in \{1,\ldots,m-1\}$,
the norm
of $U_n\circ \cdots \circ U_{2}\circ \widehat{U}_{1}$ is 
equal to $\sqrt{\widehat{a}_{1,n}}$ with
\begin{equation}\label{e:defain-bias2}
\widehat{a}_{1,n} = 
\begin{cases}
\sup_p \widetilde{\beta}_{1,n,p}^2 + \eta_{1,n}^2  & \mbox{in case \ref{caseNN1}}\\
\sup_p \left((\beta_{1,n,p}+\widetilde{\beta}_{1,n,p})^2\right) + \eta_{1,n}^2  & \mbox{in case \ref{caseNN2}}
\end{cases}
\end{equation}
and the norm of $\overline{U}_m\circ \cdots \circ U_{2}\circ \widehat{U}_{1}$
is  equal to $\sqrt{\widehat{a}_{1,m}}$ with
\begin{equation}\label{e:defain-projm2}
\widehat{a}_{1,m} = 
\begin{cases}
\sup_p
\widetilde{\beta}_{1,m,p}^2 & \mbox{in case \ref{caseNN1}}\\
\sup_p
(\beta_{1,m,p}+ \widetilde{\beta}_{1,m,p})^{2}  & \mbox{in case \ref{caseNN2}.}
\end{cases}
\end{equation}
\end{lemme}

\begin{proposition}
\label{p:LipVNN-proj2}
Let $m \in \mathbb{N}\setminus \{0,1\}$. 
For every $i\in \{2,\ldots,n\}$ and $n\in \{1,\ldots,m-1\}$,
let $a_{i,n}$ be defined by
\eqref{e:defain-bias} and
let $\overline{a}_{i,m}$ be given by
\eqref{e:defain-proj}.
For every $n\in \{1,\ldots,m\}$, let $\widehat{a}_{1,n}$
be defined by \eqref{e:defain-bias2} and  \eqref{e:defain-projm2}.
Define $(\widehat{\theta}_n)_{1\le n \le m}$ recursively by
\begin{align}
&(\forall n \in \{1,\ldots,m-1\})
\quad 
\widehat{\theta}_n = \sqrt{\widehat{a}_{1,n}}+
\sum_{i=2}^n \widehat{\theta}_{i-1} \sqrt{a_{i,n}}\;,\\
&\widehat{\theta}_m = \sqrt{\widehat{a}_{1,m}}+
\sum_{i=2}^m \widehat{\theta}_{i-1} \sqrt{\overline{a}_{i,m}}\;.
\end{align}
Then
network \eqref{e:NNreal} is $\widehat{\theta}_m /2^{m-1}$-Lipschitz.
\end{proposition}

The averagedness properties of network \eqref{e:NNreal} in cases  \ref{caseNN1}
and \ref{caseNN2} are consequences of these results.

\begin{proposition}
\label{prop:ab-proj}
Let $m \in \mathbb{N}\setminus \{0,1\}$. 
Let $\widehat{a}_{1,m}$ be defined in Lemma \ref{e:normUni-proj2} and $\overline{\theta}_m$ be defined in Proposition \ref{p:LipVNN-proj2}.
Let $\alpha \in [1/2,1]$.
Define
\begin{equation}
\widehat{b}_\alpha =   
\begin{cases}
\sup_p
(\widetilde{\beta}_{1,m,p}-2^m(1-\alpha))^2 & \mbox{in case \ref{caseNN1}}\\
\sup_p
(\beta_{1,m,p}+ \widetilde{\beta}_{1,m,p}-2^m(1-\alpha))^{2}  & \mbox{in case \ref{caseNN2}.}
\end{cases}
\end{equation}
 If
 \begin{equation}\label{e:condalphanew-proj2}
\sqrt{\widehat{b}_\alpha} -\sqrt{\widehat{a}_{1,m}} \leq 2^m \alpha - 2\widehat{\theta}_m\;,
\end{equation}
then 
network \eqref{e:NNreal} is $\alpha$-averaged.
\end{proposition}

\begin{preuve}
Let us calculate
the operator norms of $ \widehat{U} =  \overline{U}_m \circ U_{m-1}\cdots
\circ U_{2}\circ \widehat{U}_{1}.$ and 
$ \widehat{U} - 2^m (1-\alpha) \1 $.
Applying Lemma~\ref{e:normUni-proj2} when
$n=m$ yields
\[
\| \widehat{U} \|^2 = \widehat{a}_{1,m} 
\; .
\]
By following the same reasoning as in the proof of Proposition \ref{prop:ab-bias}, 
we get
\begin{equation}
\|\widehat{U} - 2^m (1-\alpha)\1\|^2
= \widehat{b}_{\alpha}.
\end{equation}
By applying now Proposition~\ref{prop:oldres}\ref{prop:oldresii},
we deduce that, when Condition \eqref{e:condalphanew-proj2} holds, 
network~\eqref{e:NNreal} is $\alpha$-averaged.
\end{preuve}

%%%%%%%%%%%%%%%%%%%%%%%%%%%%%%%%%%%%%%%%%%%%%%%%%%%%%%%%%%%%%%%%%%%%%%%%%%%%%%%%%%%%%%%%%%%%%%%%%%%%%%%%%%%%%%%%%%%%%%%%%%%%%%%%%
%%% SECTION 3
\section{Numerical Examples}
\label{section:numerics}
%%%%%%%%%%%%%%%%%%%%%%%%%%%%%%%%%%%%%%%%%%%%%%%%%%%%%%%%%%%%%%%%%%%%%%%%%%%%%%%%%%%%%%%%%%%%%%%%%%%%%%%%%%%%%%%%%%%%%%%%%%%%%%%%%

In this section, we present numerical tests 
carried out in the case of the class of Abel integral operators.
We present in more details the architecture chosen to build the neural network. Several numerical examples are provided to illustrate the accuracy of the proposed method.
The stability of the neural network is evaluated by computing its Lipschitz constant by relying upon the results of Section~\ref{section:theory}.

%
%%%%%%%%%%%%%%%%%%%%%%%%%%%%%%%%%%%%%%%%%%
\subsection{Problem formulation}
To implement the neural network defined by~\eqref{def:modelNN}, we focus on the Abel integral operator
\begin{equation}
	   \label{def:Tabel}
  	\begin{array}{cc|ccc}
  	 T & : & L^2(0,1)  & \to & L^2(0,1) \\
  	   &   & x   & \to & 
	 \displaystyle y(t) = \frac{1}{\Gamma(a)}\int_0^{t} (t-s)^{(a-1)} x(s) \d s \; ,\\
  	\end{array}
 \end{equation}
where $a>0$ and $\Gamma$ is the classical Gamma function,
$\Gamma(a)= \int_0^{+\infty} t^{a-1} e^{-t} \d t$. 
The Abel operator $T$ is injective, linear,  and compact.
The inverse problem linked to the Abel transform has been widely studied from a theoretical viewpoint, as in~\cite{gorenflo1997fractional}.
The range of $T$ is a subset of $H^{-a}(0,1)$, the dual space of~$H^a(0,1)$, and the problem is ill-posed of order $a$ in the sense of~\cite{Engl1996}.
 
 Recovering $x$ from a noisy measurement $\yd= Tx+v^\delta$ is an inverse problem linked to a large variety of experimental contexts in physics.
Indeed, 
the operator $T$ allows to define derivatives of fractional order for $a<1$ and integrals of arbitrary order for $a>1$.
The most common case
 is the semi-derivative, when $a=1/2$.
Typically, the inverse problem consists in searching a distribution
of a two-dimensional or three dimensional object
from measurements of the projection of this quantity onto an axis,
in which case the radial distribution is linked to the values of these projections via the Abel transform
(see plasmas and flames~\cite{aakesson2008}, tomography~\cite{Dribinski2002}, or astrophysics~\cite{kumar2015analytical}
).
In a different context, fractional calculus appears to be very convenient to describe properties of polymers~\cite{podlubny1998} or surface-volume reaction problems~\cite{evans2017applications}.
Subsequently, a large number of physical applications have been documented in~\cite{gorenflo1997fractional}.

According to the theory in~\cite{Gorenflo1991}, the derivative operator is given as a power of the Laplacian denoted by $B$, defined on $\mathcal{D}(B)$:
\begin{equation}
	   \label{def:D}
  \left\{
\begin{array}{rl}
B = & \; - \Delta \\ 
 \mathcal{D} (B) = & \; \left\{ x\in H^{2 }(0,1) \mid 
               x (1) =0 \;, \; \;
               x' (0) =0 \;  \right\}
               \; . 
\end{array}
\right.
 \end{equation}
Then,  
the continuous derivative operator $D$ in~\eqref{def:varJ} is chosen as $D = B^{r/2}$, with $r>0$ characterizing the order of derivation.
This choice ensures that the continuous operators $T^*T$ and $D^*D$ commute,  since for $x\in L^2(0,1)$, we have $BT^*Tx =x$ and for $x \in \mathcal{D}(B)$, we have $T^*TBx = x$.
%
%%%%%%%%%%%%%%%%%%%%%%%%%%%%%%%%%%%%%%%%%%
\subsection{Discretization}\label{se:discretAbel}

We first describe the discretization choices to pass from our continuous framework to a numerical setting.
Network~\eqref{def:modelNN} is made up of continuous operators. 
To carry out our experiments, we propose the following discretization.
We suppose that the measured signal $y = Tx$ is acquired on a regular mesh of $N$ points, $(t_i)_{0  \leq i \leq N-1}$ in the interval $[0,1]$, 
with $t_0=0$ and $t_{N} =1$. 
The measured signal $ y =  (y_i)_{0  \leq i \leq N-1} =(y(t_i))_{0  \leq i \leq N-1}$ belongs to the space endowed with the finite element basis 
$(e_i)_{0\leq i\leq N-1}$ 
associated to $(t_i)_{0  \leq i \leq N-1}$. 
However, instead of working only in the finite element basis, 
we also consider projection of the signal onto the span of the first $K$ eigenvectors of the self-adjoint nonnegative operator $T^*T$.
This choice is justified for two reasons. 
First, in such basis, the discretized forms of
operators $T^* T $ and $ D^* D $ respectively defined by~\eqref{def:Tabel} and~\eqref{def:D} are diagonal and therefore commute. 
This is a prerequisite to apply Propositions~\ref{p:LipVNN-proj2} and~\ref{prop:ab-proj}. 
Second, the eigenvectors of $T^*T$ 
denoted by $(u_k)_{k \in \N}$ 
are trigonometric polynomials 
and the retained discretization method is a spectral one, as defined in~\cite{Gottlieb1977}, or~\cite{Boyd2001}.
This discretization method can fully account for the regularity of the initial condition on $x$,
under extra mild assumptions.

We denote by $(u_k, \beta_{T,k})_{k \in \N} $ the eigensystem of $T^*T$. 
Note that since $T$ is a compact operator, the eigenvectors $(u_k)_{k \in \N}$ is a set of orthonormal eigenvectors, and $(\beta_{T,p})_{n \in \N}$ are strictly positive eigenvalues.
The signal $y$ is then discretized in the basis formed by the $K$ first eigenvectors $(u_k)_{0\leq k\leq K-1}$.
Explicit values of $(u_k, \beta_{T,k})$ are given in~\cite{Gorenflo1999}.
As already stated, the operators $T^*T$ and $D^*D$ reduce to diagonal matrices, 
with the following eigenvalues on their diagonal:  
\begin{equation}
(\forall k \in \{0,\ldots,K-1\})\quad
    \beta_{T,k} = \left( \frac{4}{\pi^2(1+2 k)^2}\right)^{a}
    \; ,
    \quad
    \beta_{D,k} = \left( \frac{\pi^2(1+2 k)^2}{4} \right)^r
    = \beta_{T,k}^{-r/a}
    \; .
\end{equation}
Hereafter, we consider that $r=1$, and $D = B^{1/2}$, with $B$ defined by~\eqref{def:D}.
We compute the change basis matrix denoted by $P=(P_{i,j})_{0\le i,j\le N-1}$ with, for 
every $i\in \{0,\ldots,N-1\}$ and $j\in \{0,\ldots,N-1\}$,
\begin{equation}
P_{i,j} = 
\frac{2 \sqrt{2}}{\gamma_j} \text{cos}\left(\frac{2i+1}{2N}\gamma_j\right) \text{sin}\left(\frac{1}{2N}\gamma_j\right)
\; ,
\quad
\gamma_j = \sqrt{\beta_{T,i}} 
         = \left( \frac{2}{\pi(1+2 n)}\right)^{a}
\; .
\end{equation}
The operator $T$ does not intervene in the neural network~\eqref{def:modelNN}, as only $T^*T$ and $D^*D$ do. 
However, to generate synthetic data and the associated bias $b_0$, we need also a discretization for the operator $T$.
Therefore, $T$ is approximated by $T_{\text{elt}}$ as a computation of an integral using the trapezoïdal rule, with stepsize $h=1/(N-1)$, and, for $a \neq 1$,
for $0\leq i<N$, $0\leq j<N$,
\begin{equation}
    \label{def:Telt}
(T_{\text{elt}})_{i,j} = 
\begin{cases}
\; \displaystyle \frac{1}{\Gamma(a)} \frac{h^a}{2a}
     \left( (i-j+1)^{a} - (i-j-1)^{a}  \right) & \mbox{if $j<i$} \; ,\\
     \\
     \; \displaystyle \frac{1}{\Gamma(a)}  \frac{h^a}{2a}  (i^a - (i-1)^a) & \mbox{if $j=0, \; i \neq 0$}   \; ,\\
\\
\; \displaystyle \frac{1}{\Gamma(a)}  \frac{h^a}{2a}   & j=i , \; \mbox{if $i \neq 0$} \; ,\\
\\
\; 0 &\mbox{if $ i=j=0$, or $j > i$}.
\end{cases}
\end{equation}
Then, the operators $T$ and $T^*$ in the eigen basis are respectively approximated by $T_{\text{eig}} = P T_{\text{elt}} $ and $(T^*)_{\text{eig}} = P T_{\text{elt}}^\top $.
Thus, 
on the one hand the synthetic data are calculated and stored in the basis of the finite elements,
and, on the other hand, the algorithm operates in the basis of eigenvectors,
except for the proximity operator, 
for which a change of basis is performed before and after.

To carry out the numerical experiments, we set $N= 2 \times 10^3$ and $K=50$. 
Therefore, the regular signals $x$ are approximated by their projection onto the space generated by the first $K$ eigenvectors  of $T^*T$.

%%%%%%%%%%%%%%%%%%%%%%%%%%%%%%%%%%%%%%%%%%%%%%%%%%%%%%%%%%%%%%
\subsection{Neural network architectures and characteristics}

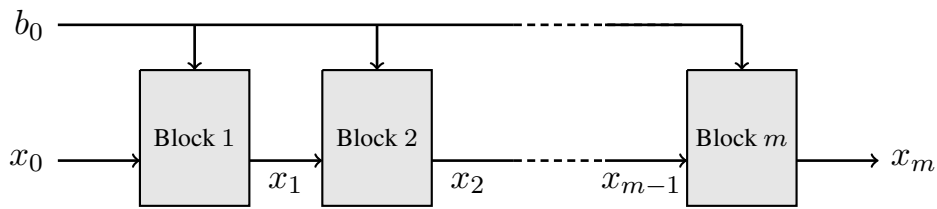
\begin{figure}[h]
\centering
\begin{tikzpicture}[scale=1.2, every node/.style={scale=1.2}]
% Style
\tikzstyle{block}=[thick, fill=gray!20]
% Initialisation
\draw [line width=1pt] (0,1.5) node[left]{$b_0$} ; 
\draw [line width=1pt] (0,0) node[left]{$x_0$} ; 
% Fleche du biais
\draw [line width=1pt] (0,1.5) -- (5,1.5) ;
% ========================================
% BLOCK #1
\draw [line width=1pt][->](0,0) -- (0.9,0) ; % x input block
\draw [line width=1pt][->](1.5,1.5) -- (1.5,1) ;   % biais
\draw [block] (0.9,1) -- (2.1,1) -- (2.1,-0.5) -- (0.9,-0.5) -- (0.9,1) ;
\node[font=\scriptsize] (B1) at (1.5,0.25) {Block $1$};
\draw [line width=1pt][->](2.1,0) -- (2.9,0) ; % x output
\draw (2.5,0) node[below]{$x_1$} ; 
% ========================================
% BLOCK #2
\draw [line width=1pt][->](3.5,1.5) -- (3.5,1) ;   % biais
% Block
\draw [block] (2.9,1) -- (4.1,1) -- (4.1,-0.5) -- (2.9,-0.5) -- (2.9,1) ;
\node[font=\scriptsize] (B2)at(3.5,0.25) {Block $2$};
\draw [line width=1pt](4.1,0) -- (5,0) ; % x output
\draw (4.5,0) node[below]{$x_2$} ; 
% ========================================
% BLOCK (...)
\draw [line width=1pt,densely dashed](5,0) -- (6,0) ; %x (...)
\draw [line width=1pt][->](6,0) -- (6.9,0) ; % x fleche
\draw (6.4,0) node[below]{$x_{m-1}$} ; 
\draw [line width=1pt, densely dashed](5,1.5) -- (6.9,1.5) ; % biais (...)
% Block
\draw [block] (6.9,1) -- (8.1,1) -- (8.1,-0.5) -- (6.9,-0.5) -- (6.9,1) ;
\node[font=\scriptsize] (B3)at(7.5,0.25) {Block $m$};
\draw [line width=1pt] (6,1.5) -- (7.5,1.5) ; % bias till end
\draw [line width=1pt][->](7.5,1.5) -- (7.5,1) ;% biais fleche
% ========================================
% END BLOCK 
\draw [line width=1pt][->](8.1,0) -- (9,0) ; % x fleche
\draw (9,0) node[right]{$x_{m}$} ; 
\end{tikzpicture}
\caption{Global architecture of neural network~\eqref{def:modelNN}}
\label{fig:structure}
\end{figure}

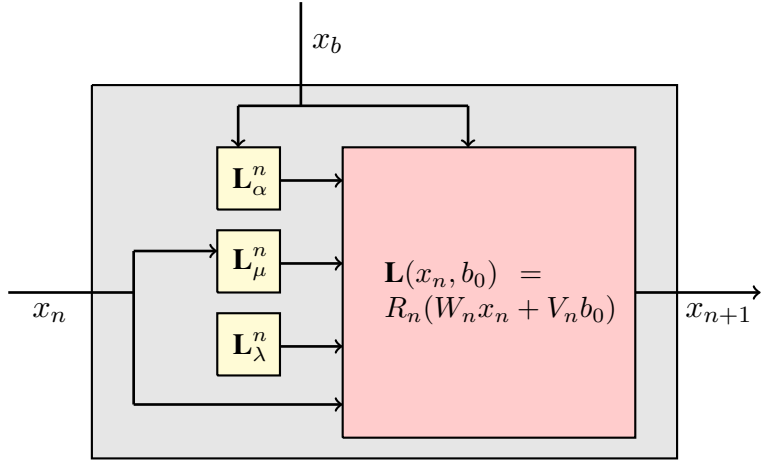
\begin{figure}
    \centering
\begin{tikzpicture}[scale=1.1, every node/.style={scale=1.1}]
% Style
\tikzstyle{block}=[thick, fill=gray!20]
\tikzstyle{iterates}=[thick, fill=red!20]
\tikzstyle{hyper}=[thick, fill=yellow!20]
% ========================================
%
\draw [block] (0,5) -- (7,5) -- (7,0.5) -- (0,0.5) -- (0,5) ;
% ========================================
%
\draw [line width=1pt](-1,2.5) -- (0.5,2.5) ; % x input block
\draw [line width=1pt] (0.5,1.15) -- (0.5,3) ; % x vertical line
\draw (-0.5,2.5) node[below]{$x_{n}$} ; 
\draw [line width=1pt](2.5,6) -- (2.5,4.75) ; % b input block
\draw [line width=1pt](1.75,4.75) -- (4.5,4.75) ; % b horizontal line
\draw (2.5,5.5) node[right]{$x_{b}$} ; 
% ========================================
% Hyper Parameter
 % gradient step (bottom)
\draw[hyper] (1.5,1.5) -- (1.5,2.25) -- (2.25,2.25) -- (2.25,1.5) -- (1.5,1.5) ;
\node[font=\small] (lambda)at(1.9,1.85) {$\textbf{L}_\lambda^n$};
\draw [line width=1pt][->] (2.25,1.85) -- (3,1.85) ; % to principal layer
% barrier (middle)
\draw[hyper] (1.5,2.5) -- (1.5,3.25) -- (2.25,3.25) -- (2.25,2.5) -- (1.5,2.5) ;
\node[font=\small] (mu)at(1.9,2.85) {$\textbf{L}_\mu^n$};
\draw [line width=1pt][->] (0.5,3) -- (1.5,3) ; 
\draw [line width=1pt][->] (2.25,2.85) -- (3,2.85) ; % to principal layer
% Regularisation (top)
\draw[hyper] (1.5,3.5) -- (1.5,4.25) -- (2.25,4.25) -- (2.25,3.5) -- (1.5,3.5) ;
\node[font=\small] (alpha)at(1.9,3.85) {$\textbf{L}_\alpha^n $};
\draw [line width=1pt][->] (1.75,4.75) -- (1.75,4.25) ; % b to barrier layer
\draw [line width=1pt][->] (2.25,3.85) -- (3,3.85) ; % to principal layer
% ========================================
% Iteration Layer
\draw [line width=1pt][->] (0.5,1.15) -- (3,1.15) ; % x to tayer
\draw [line width=1pt][->] (4.5,4.75) -- (4.5,4.25) ; % b to layer
\draw[iterates] (3,0.75) -- (3,4.25) -- (6.5,4.25) -- (6.5,0.75) -- (3,0.75) ;
\node[font=\small,text width=3cm] (alpha)at(5,2.5) {$ \textbf{L} (x_n,b_0)= R_n (W_n x_n + V_n b_0)$};
% ========================================
% Output
\draw[line width=1pt][->] (6.5,2.5) -- (8,2.5) ;
\draw (7.5,2.5) node[below]{$x_{n+1}$} ; 
\end{tikzpicture}
\caption{Architecture of one iteration - Block $n$.}
\label{fig:block}
\end{figure}

\paragraph{Structure}
The architecture that we propose here reflects the proposed Model~\eqref{def:modelNN}, that unfolds the forward-backward algorithm for minimizing functional $J$ defined by~\eqref{def:varJ} 
over a finite number of iterations $m$. 
The main difference with the classical forward-backward algorithm 
lies in the fact that only a finite number of iterations $m$ is performed 
instead of pursuing the iterations until convergence. 
Then, $m$ corresponds to the number of layers of the neural network. 
Similarly to the work proposed in~\cite{Corbineau2020}, 
each layer of the neural network consists of a block made up of hidden layers which calculate the hyper-parameters and an iteration of the forward-backward algorithm.
Here, the bias $b_0$ is taken as the discretization of $T^*\yd$, namely $PT_{\text{elt}}^\top \yd $.
% We choose to set the input $x_0$ as $B T^* \yd $, since $BT^*T = \1$. The noise in the signal $\yd$ is then amplified by this computation, since the problem is ill-posed.
The hyper-parameters are defined independently across the network in order to provide more flexibility.
The overall structure of the network is shown in Figure~\ref{fig:structure}.
The structure of its hidden layers is detailed in Figure~\ref{fig:block}.
% ========================================================================
% ==========================================================================
%% GLOBAL STRUCTURE
% ==========================================================================
% ==========================================================================

% ==========================================================================
%% ONE BLOCK
% ==========================================================================
% ========================================================================
The activation function, namely operator $R_n$ in Figure~\ref{fig:block}, 
corresponds the proximal operator associated with~$g$ appearing  in~\eqref{def:varJ}.
We remind that, for an indicator function of a nonempty closed convex set, the proximity operator is a projection onto this set.
However, such activation functions, 
especially in the case where one wishes to guarantee the positivity of the solution, 
may show bad properties during gradient back-propagation and training.
These include vanishing gradient problems
as shown in~\cite{misra2019},~\cite{Pedamonti2018}
for the Rectified Linear Unit (ReLU) function.
Then, we choose to consider instead a logarithmic barrier $g$
to enable prior knowledge in the algorithm,
as proposed in~\cite{Corbineau2020}.
The activation is no more constant and depends on the gradient step $\lambda_n$ and the barrier parameter $\mu_n >0$ as 
\[
R_n = \text{prox}_{\lambda_n\mu_n g} \; .
\]

\paragraph{Constraint and proximity operator}
More precisely, 
we experiment two possible choices for
the function $g$ in expression~\eqref{def:varJ}.
%The function is not chosen as the indicator function of a bounded convex set $C$.
%Instead, we enable 
As mentioned above, the prior knowledge on the constraint set $C$ is thus embedded in the network through the logarithmic barrier function $g$.
In both cases, 
\begin{equation}
    \label{def:logbarrier}
    \begin{cases}
    C = \{ x \in L^2(0,1) \mid c_i(x) \geq 0, \; 1 \leq i \leq p \} \; ,\\
    \\
    \text{int} \; C = \{ x \in L^2(0,1) \mid \; c_i(x) > 0, \; 1 \leq i \leq p \} \; ,\\
    \\
    (\forall x \in L^2(0,1))\quad
    g (x) = \left\{
    \begin{array}{cc}
       - \sum_{i=1}^p  \text{ln} \; (c_i(x)) & \text{if} \; x \in \text{int} \; C \\
        +\infty                              &  \text{otherwise} \; ,
    \end{array}
    \right.
    \end{cases}
\end{equation}
where $(c_{i})_{1\le i \le p}$ are suitable functions allowing us to describe the constraint set.
First, we consider that the signal $x$ has a minimum value $x_{\text{\rm min}}$ and a maximum value $x_{\text{\rm max}}$.
Then $C$ can be rewritten as 
\begin{equation}
\label{constr:cube}
C = \left\{ 
x \in L^2(0,1) \mid x\geq x_{\text{\rm min}}, \; -x \geq -x_{\text{\rm max}}
\right\} \; .
\end{equation}
This kind of constraint can be useful for example when the signal
the experimenter wishes to recover 
corresponds to a positive and bounded physical quantity. 
Secondly, we consider an affine constraint such as, for $j>0$,
\begin{equation}
\label{constr:slab}
C = \left\{ 
x \in L^2 \mid  \;  0 \leq \int_0^1 t^j x(t) \d t \leq 1
\right\} \; . 
\end{equation}
This constraint reflects the fact that a physical quantity linked to the signal is bounded. 
For $j\in \{1,2,3\}$, the moment of order $j$ involved in \eqref{constr:slab} represents the total mass of elements in a 1D, 2D, or 3D system, respectively.

\begin{figure}[h]
\centering
\begin{tikzpicture}[scale=0.8, every node/.style={scale=1.2}]
% Input 
\draw [font=\scriptsize,text width=1.5cm] (-0.5,3) node[right]{$x_n$} ; 
\draw [font=\tiny,text width=1.5cm] (0.5,6) node[above]{size $(1\times n)$}; ; 
\draw (0.25,0) grid[step=0.25] (0.5,6) ;% grille 1
% Convolution 1
\draw [font=\scriptsize,text width=1.5cm] (1,-0.25) node[right]{AvgPool +Softplus} ; 
\draw (2.99,0.99) grid[step=0.25] (3.25,5) ; % grille 2
\draw [font=\tiny,text width=1.5cm] (3.5,5) node[above]{$(1\times 64)$};
\draw [dotted,thick] (0.5,6) -- (3,5); %top
\draw [dotted,thick] (0.5,0) -- (3,1); %bottom
%%%
% Kernel
\draw[red] (0.25,3.99) grid[step=0.25] (0.5,5) ;
\draw [dotted,thick] (0.5,5) -- (3,4.5); %top
\draw [dotted,thick] (0.5,4) -- (3,4.25); %bottom
%%%
% Convolution 2
\draw  [font=\scriptsize,text width=1.5cm] (3.5,0.5) node[right]{AvgPool +Softplus} ; 
\draw (5.49,1.49) grid[step=0.25] (5.75,4.5) ;
\draw [font=\tiny,text width=1.5cm] (6,4.5) node[above]{$(1\times 16)$};
\draw [dotted,thick] (3.25,5) -- (5.5,4.5); %top
\draw [dotted,thick] (3.25,1) -- (5.5,1.5); %bottom
% Fully connected
\tikzstyle{neuron}=[fill=gray!20]
\node (n1) at (7.5,4) {};
\node (n2) at (7.5,3.5) {};
\node (n3) at (7.5,2.5) {};
\node (n4) at (7.5,2) {};
\draw [neuron] (n1) circle(0.2);
\draw [neuron] (n2) circle(0.2);
\draw [neuron] (n3) circle(0.2);
\draw [neuron] (n4) circle(0.2);
\node (x1) at (5.75,4) {};
\node (x2) at (5.75,2) {};
\draw[dotted] (x1) -- (n1) ;
\draw[dotted] (x1) -- (n2) ;
\draw[dotted] (x1) -- (n3) ;
\draw[dotted] (x1) -- (n4) ;
\draw[dotted] (x2) -- (n1) ;
\draw[dotted] (x2) -- (n2) ;
\draw[dotted] (x2) -- (n3) ;
\draw[dotted] (x2) -- (n4) ;
\node (nout) at (8.5,3) {};
\draw [neuron] (nout) circle(0.2);
\draw[dotted] (n1) -- (nout) ;
\draw[dotted] (n2) -- (nout) ;
\draw[dotted] (n3) -- (nout) ;
\draw[dotted] (n4) -- (nout) ;
\draw[->] (nout) -- (10,3) ;
\draw [font=\scriptsize,text width=1.5cm] (10,3) node[below]{Softplus};
\draw [font=\scriptsize,text width=1.5cm,text badly ragged] (6,1.25) node[right]{Fully connected} ; 
\end{tikzpicture}
\caption{Architecture of one hidden layer $\textbf{L}_\mu^n$ computing the barrier parameter.}
\label{fig:L_mu}
\end{figure}
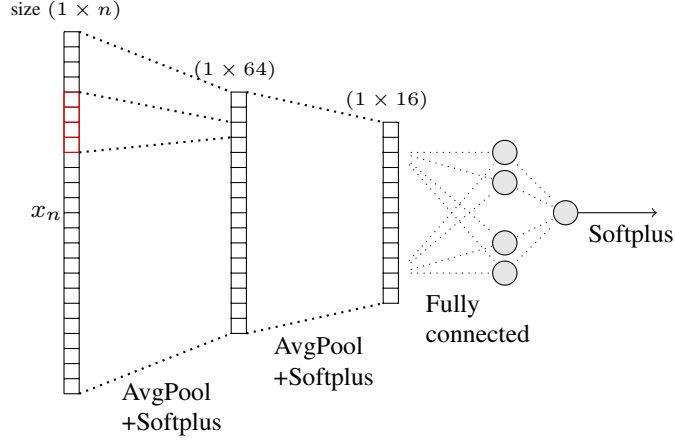

The computation of the proximity operator associated to logarithmic barrier functions~\cite{Chouzenoux2012},
after discretization,
can be found in~\cite{Corbineau2020}.
In our case,
the barrier parameter in each layer $n$, denoted by $\mu_n$ 
is estimated with a convolutional neural network,
which takes as input the output $x_{n}$ of the $(n-1)$-th layer.
The detailed architecture of $\textbf{L}_\mu^n$ is depicted in Figure~\ref{fig:L_mu}.
Since the barrier parameter is positive, we enforce this constraint by an approximation of the ReLU activation function, namely Softplus~\cite{Dugas2000}, with $\beta >0$,
\begin{equation}
\label{def:softplus}
    \text{Softplus} (x,\beta) = \frac{1}{\beta} \text{ln} \left( 1+ e^{\beta x} \right)
    \; .
\end{equation}

% ========================================================================
%
% ========================================================================
%
\paragraph{Other parameters}
We introduce between each layer a hidden layer responsible for computing the gradient step size $\lambda_n$, and the regularization parameter $\tau_n$, respectively.

The gradient descent step $(\lambda_n)_{1 \leq n \leq m}$ only depends on the structure of the network. 
This parameter is then trained without any prior knowledge.
Since its value is positive, we compute it as:
\begin{equation}
    \lambda_n = \text{Softplus} (c_n),
\end{equation}
where $c _n$ is a scalar parameter of the network learned during training.

The regularization parameters $(\tau_n)_{1 \leq n \leq m}$ in neural network~\eqref{def:modelNN} should only depend on the bias $b_0$.
Indeed, for a regularization of the generalized Tikhonov type, 
the regularization parameter theoretically depends on the regularity of the a priori and on the noise level. 
The theoretical optimal value of this parameter can be explicitly computed, as shown in~\cite{natterer1984} or~\cite{Gorenflo1991},
\begin{equation}
    \tau = c \left( \frac{\delta}{\rho}\right)^{\frac{2(a+r)}{a+q}}
    \; ,
\end{equation}
where $\rho = \| x\|_{L^{q}(0,1)} $,
and $x \in H^q(0,1)$ represents the ideal signal, 
$\delta$ is the noise level in $L^2$ norm,
$a$ is the degree of ill-posedness of the inverse problem, $r$ the level of regularization (or the order of the differential term in the regularization),
and $c$ is a constant.
Since we do not have access to the noise level, we estimate it thanks to the Fourier transform of the signal. 
We assume here that the noise
corresponds to the high frequency components of the signal.
This assumption is only used to obtain an approximate value of the error. 
Subsequently,  
the algorithm makes it possible to search the optimal value
without any assumption on the Fourier spectrum of the error. 
This is achieved by learning a constant $ d_n $ such that
\begin{equation}
    \tau_n = \text{Softplus} (d_n) \; \left(
    \frac{\|b_0 - \text{FFT}_{f_{\text{max}}}(b_0)\|_2}{\|\text{FFT}_{f_{\text{max}}}(b_0)\|_q}
    \right)^{\frac{2(a+r)}{a+q}}
    \; ,
\end{equation}
where the operator $\text{FFT}_{f_{\text{max}}}$ cuts the frequencies of the Fourier transform greater than $f_{\text{max}}$, $r$ is the order of derivation in the regularization term, and $q$ is the order of regularity of the a priori, i.e. $x\in H^q(0,1)$.
This form of regularization is theorically the best choice as long as $r \geq q/2-a$,
as shown in~\cite{natterer1984}.
Moreover, this insures that the dependence of the $\tau_n$ parameter of the network on the bias $b_0$ is of second order and can be neglected while computing the Lipschitz constant.
%
%%%%%%%%%%%%%%%%%%%%%%%%%%%%%%%%%%%%%%%%%%
\subsection{Dataset and experimental settings}
\paragraph{Synthetic Data}
To ensure the universality of the approach, 
we train the network on a wide variety of functions, 
without too strong a priori on their form or their properties.
For example, we do not want to restrict our training to Gaussian-like functions 
which would be likely to be oversimplistic models.
%only
%in the fear that the learn properties might 
%not work for other type of signal $x$.
We are therefore looking for a sufficiently rich dictionary of functions sampled over $N$ points.
To create a diverse dataset of positively distributed functions, we found convenient to use histograms of color images from a standard image dataset.
However, 
in order to properly reflect the a priori of regularity, 
the following processing is then carried out to these histograms.
The functions are first smoothed using a Savitzky-Golay filter, with filter length 21 and polynomial order 5. 
%the length of the filter window is set is 21 and the order of the polynomial used to fit the samples is 5. 
Then, to ensure that such signals are in the range of $T^*T$, the outputs of the filter are padded at $t=0$ by a constant value, and at $t=1$ by zero.
Finally, the signals are projected into the eigenvector basis described in Section \ref{se:discretAbel}.
This process ensures that 
the obtained  signal $x$ in the training set belong to $C^\infty$, as the eigenvectors do.
In particular, $x$ belongs to the space  $\mathcal{D} (B)$ defined in  \eqref{def:D},
%$\left\{ x\in H^2(0,1) \mid x(1)=0, \; x'(0)=0 \right\}$,
in which case the regularization is optimal in the sense of~\cite{Engl1996} (
see~\cite{Gorenflo1999}).

In order to synthesize noisy signals $\yd$, the discrete transformation $T_{\text{elt}}$ defined by~\eqref{def:Telt} is applied to the set of signal $x$ created as aforementioned. Then, a zero-mean white Gaussian noise with a preset standard deviation $\delta$ is added. 
\begin{figure}[h]
    \centering
    \includegraphics[width=6cm]{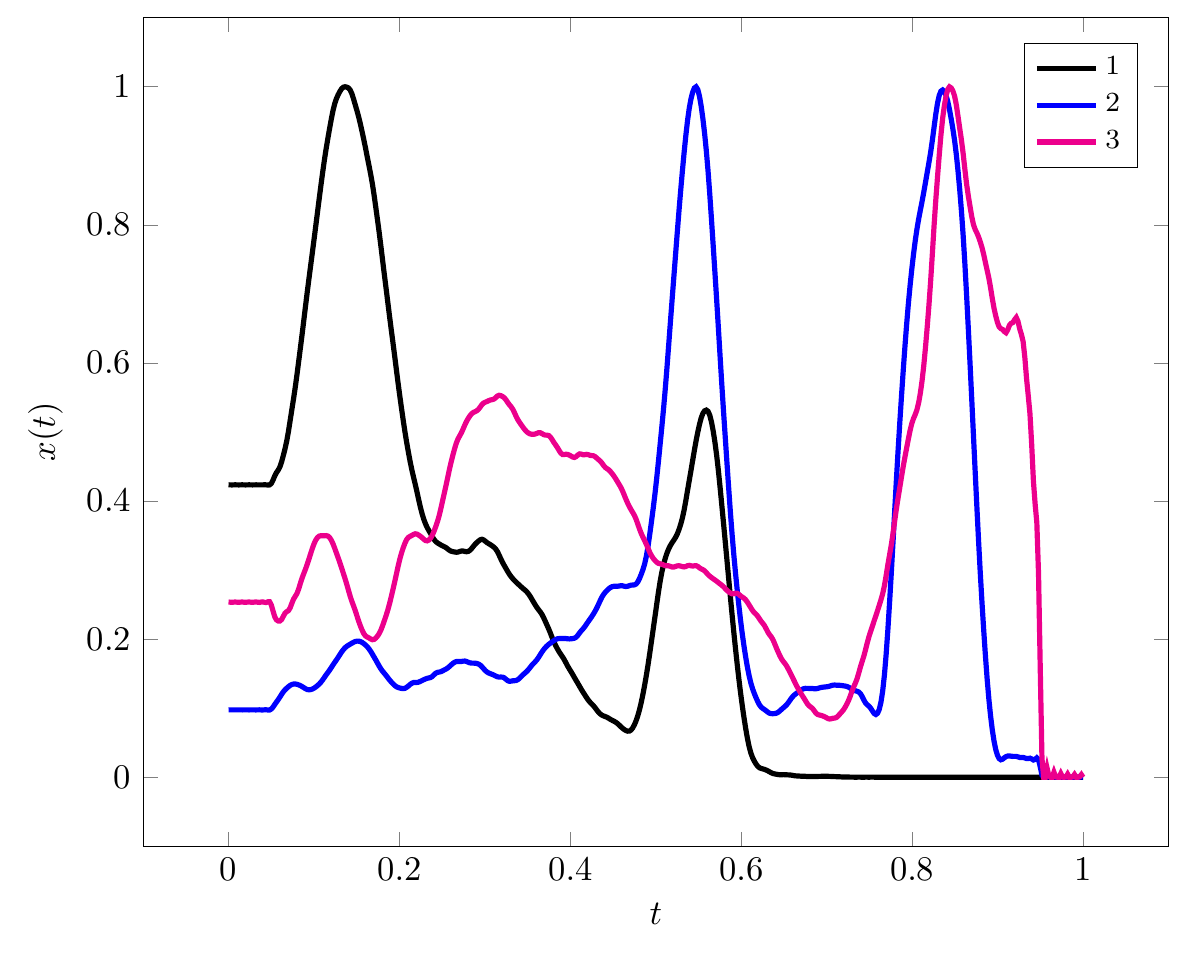}
    \caption{Example of three signals of synthetic data for the constraint~\eqref{constr:cube}. These signals are bimodal or almost bimodal, with variable peak widths. This dataset presents a great diversity of functions and demonstrates the agnostic nature of the model in order to represent a large panel of physical signals. The imposed constraints are the regularity of the signal (here $C^\infty(0,1)$) and the boundary conditions at $t=0$ and $t=1$. }
    \label{fig:dataset}
\end{figure}
\begin{figure}
    \centering
    \includegraphics[width=5.5cm,scale=0.5]{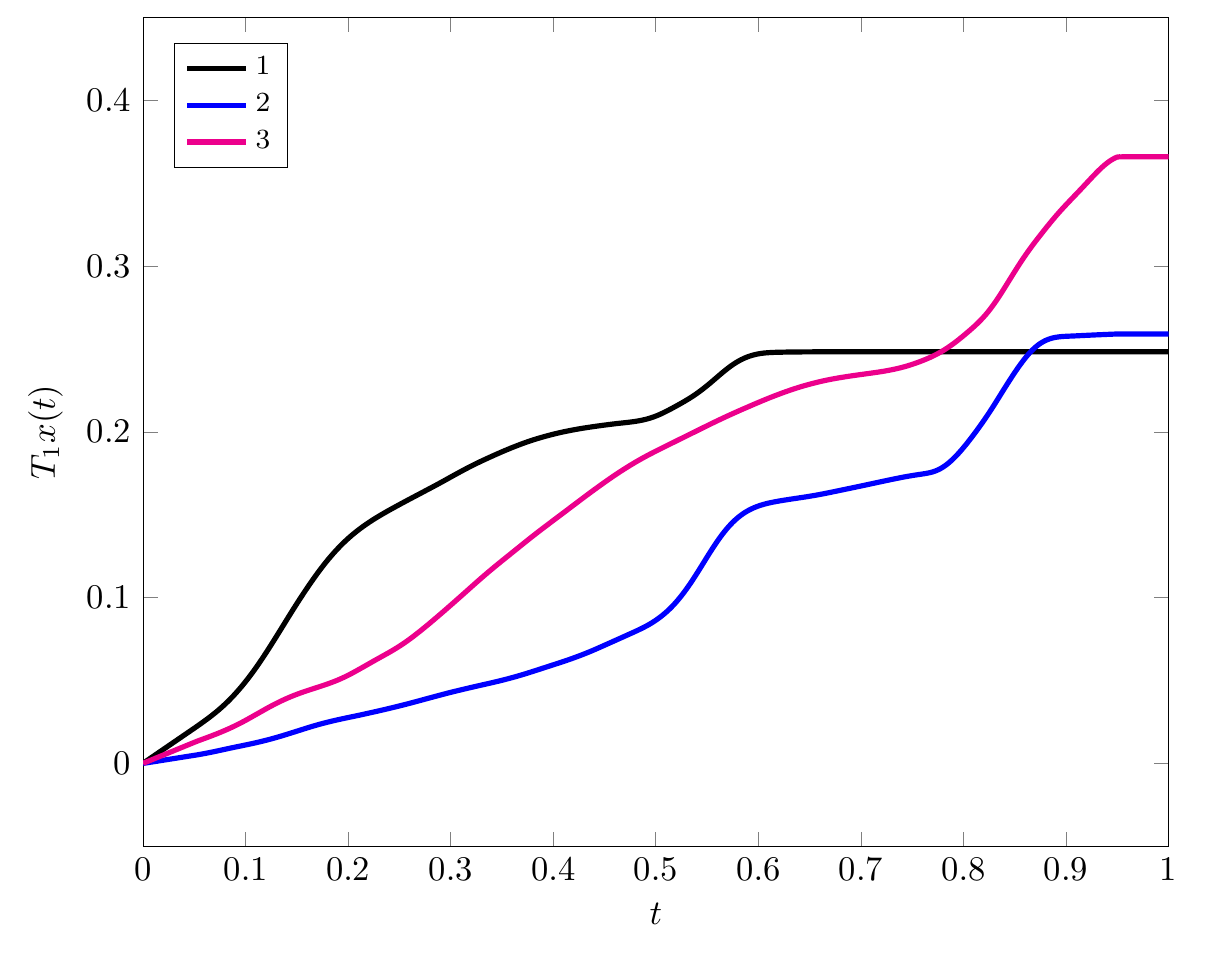}
    \includegraphics[width=5.5cm,scale=0.5]{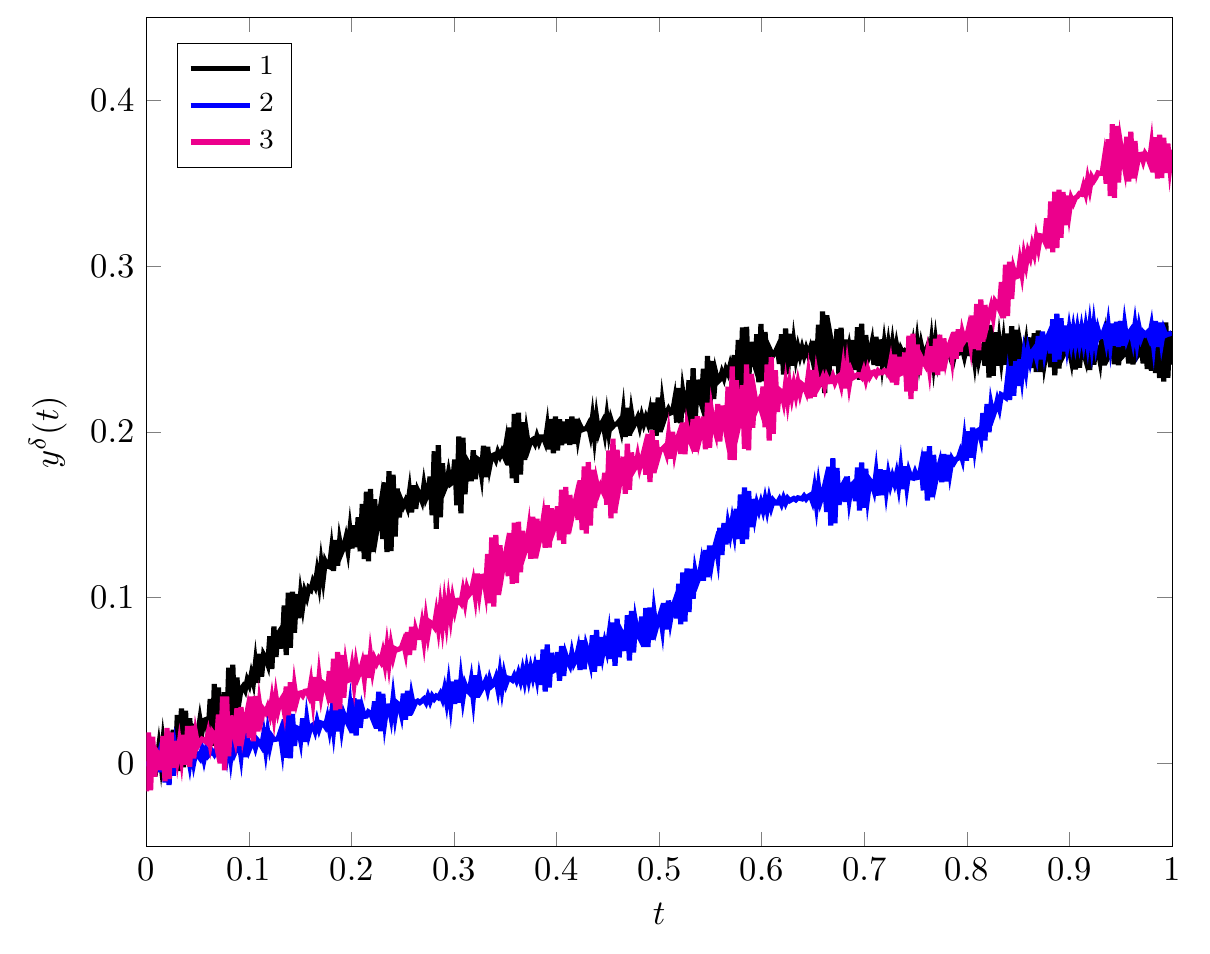}
    \caption{For the three examples displayed in Figure~\ref{fig:dataset},
    we plot on the left their image by $T$ defined by~\eqref{def:T} for $a=1$,
    and on the right the same signal after addition of noise with level $\delta = \| y^\delta - T_1 x \|_{L^2} = 0.05 \|T_1 x\|_{L^2}$.
    }
    \label{fig:yydelta-a=1}
\end{figure}
\begin{figure}
    \centering
    \includegraphics[width=5.5cm,scale=0.5]{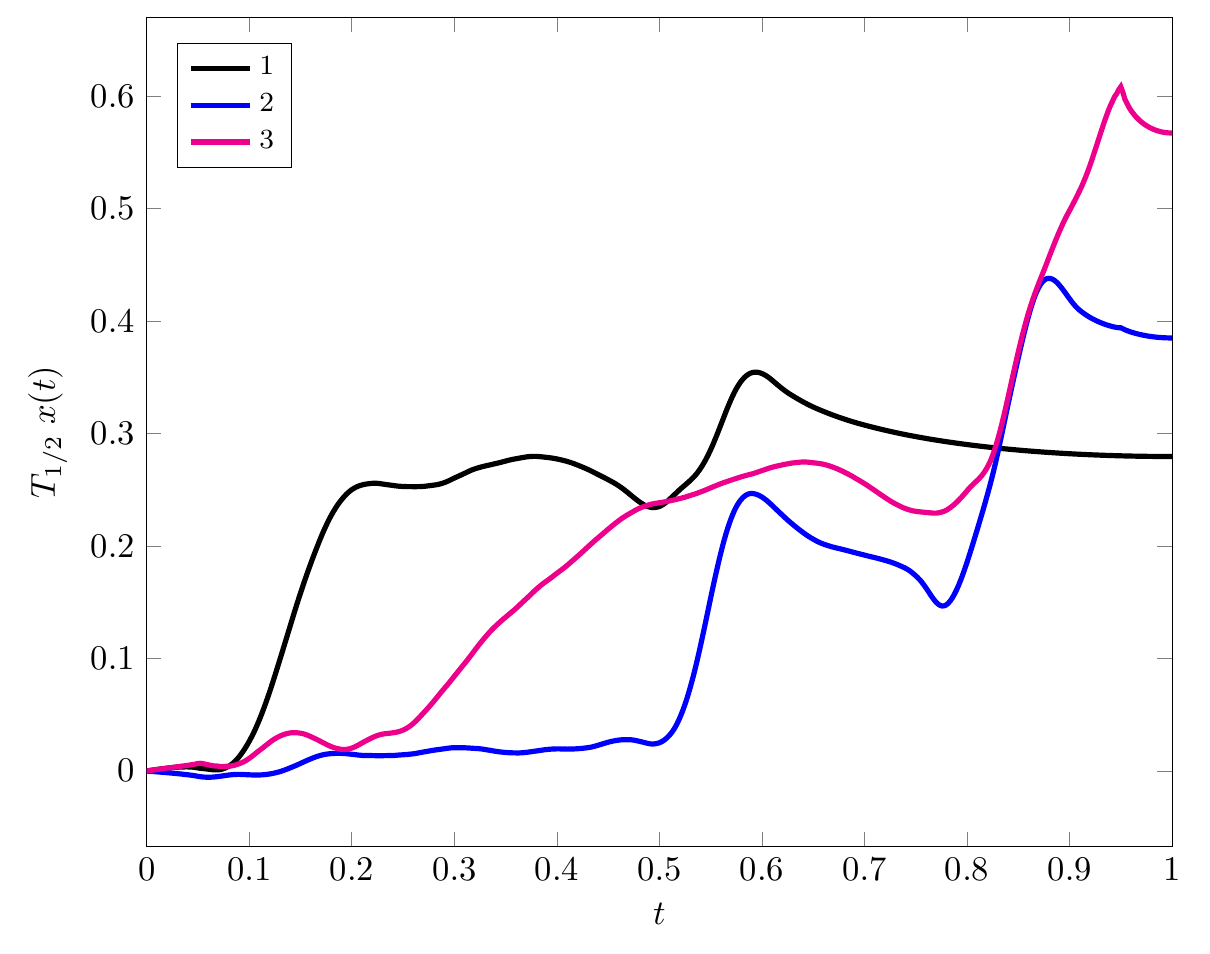}
    \includegraphics[width=5.5cm,scale=0.5]{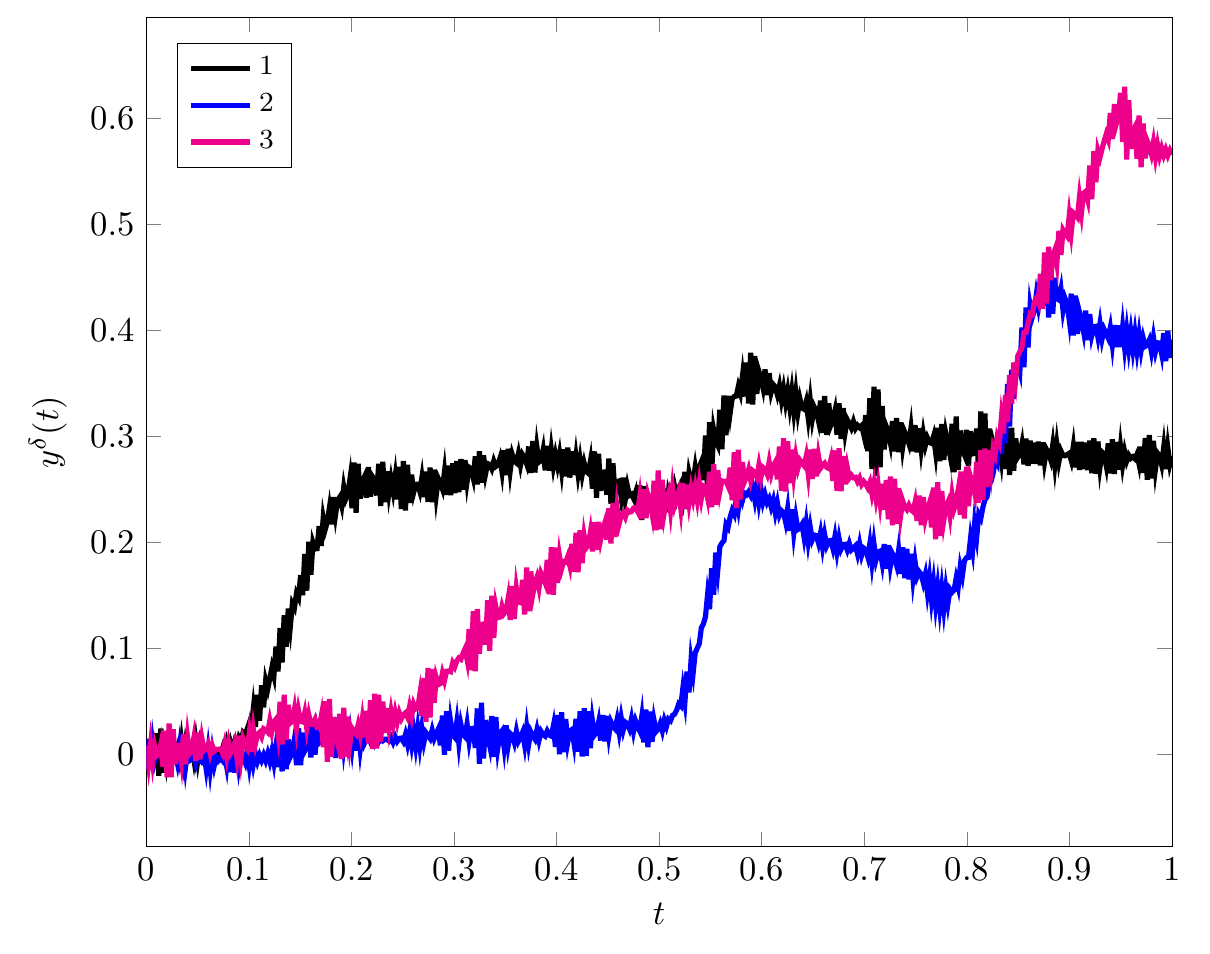}
    \caption{For the three examples displayed in Figure~\ref{fig:dataset},
    we plot on the left their image by $T$ defined by~\eqref{def:T} for $a=1/2$,
    and on the right the same signal after addition of noise with level $\delta = \| y^\delta - T_1 x \|_{L^2} = 0.05 \|T_1 x\|_{L^2}$.}
    \label{fig:yydelta-a=1/2}
\end{figure}

Figure~\ref{fig:dataset} represents an example of three signals simulated by our method.
In Figures~\ref{fig:yydelta-a=1}~and~\ref{fig:yydelta-a=1/2}, we display the image of those signals by the operator $T$ defined by~\eqref{def:T} with and without the presence of additive noise, respectively for $a=1$ and $a=1/2$.
For $a = 1$, we recover the 1D integral operator. 
All our datasets and codes implemented in Pytorch are available online~\footnote{https://github.com/ceciledellavalle/FBResNet}.

\paragraph{Training}
The network is classically trained in an end-to-end fashion. 
The gradient of the proximity operator is explicitly coded and inserted into the back-propagation according to the chain rule.
We trained the network over 30 epochs with a learning rate of $10^{-3}$, by using a training set of 400 signals. 
We use the Adam optimizer~\cite{kingma2014adam} to minimize the training loss, which is taken as the mean square error.
We compute the validation loss at every epoch by using a set of 200 signals.
The batchsize is equal to one.
The training takes approximately three to four hours on an NVIDIA Titan Xp GPU, while the computational time required for testing one signal is only about 50ms on a 2.9 GHz 6-Core Intel Core i9.

Using Proposition~\ref{p:LipVNN-proj2},
the Lipschitz constant is estimated at each epoch of training. 
One observes various behaviors, according to the initial values of the hyper-parameters. 
Either the Lipschitz constant increases until stabilizing, 
or it decreases. 

\subsection{Results and discussion}
\begin{figure}
    \centering
    \includegraphics[width=6cm,scale=0.5]{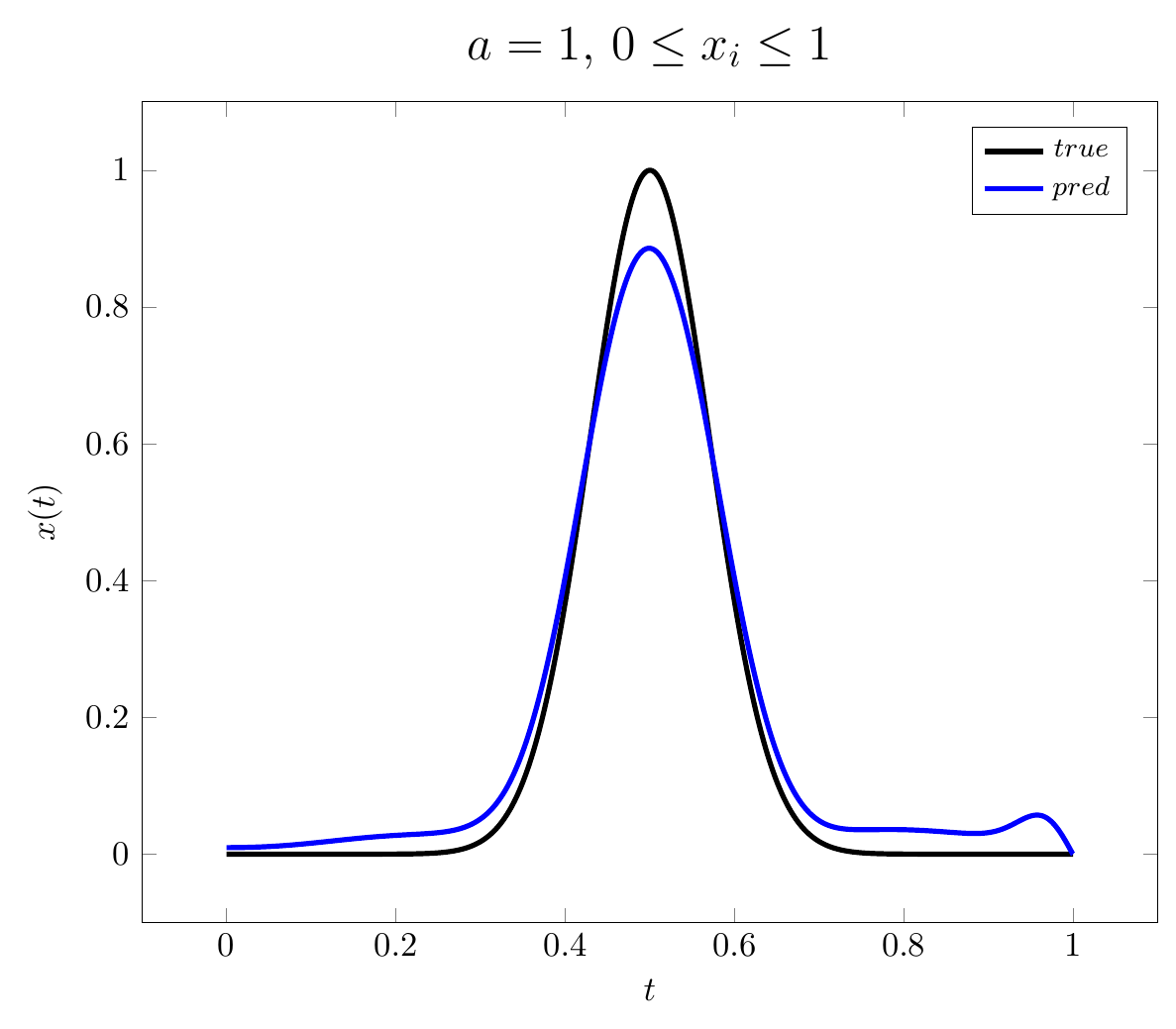}
    \includegraphics[width=6cm,scale=0.5]{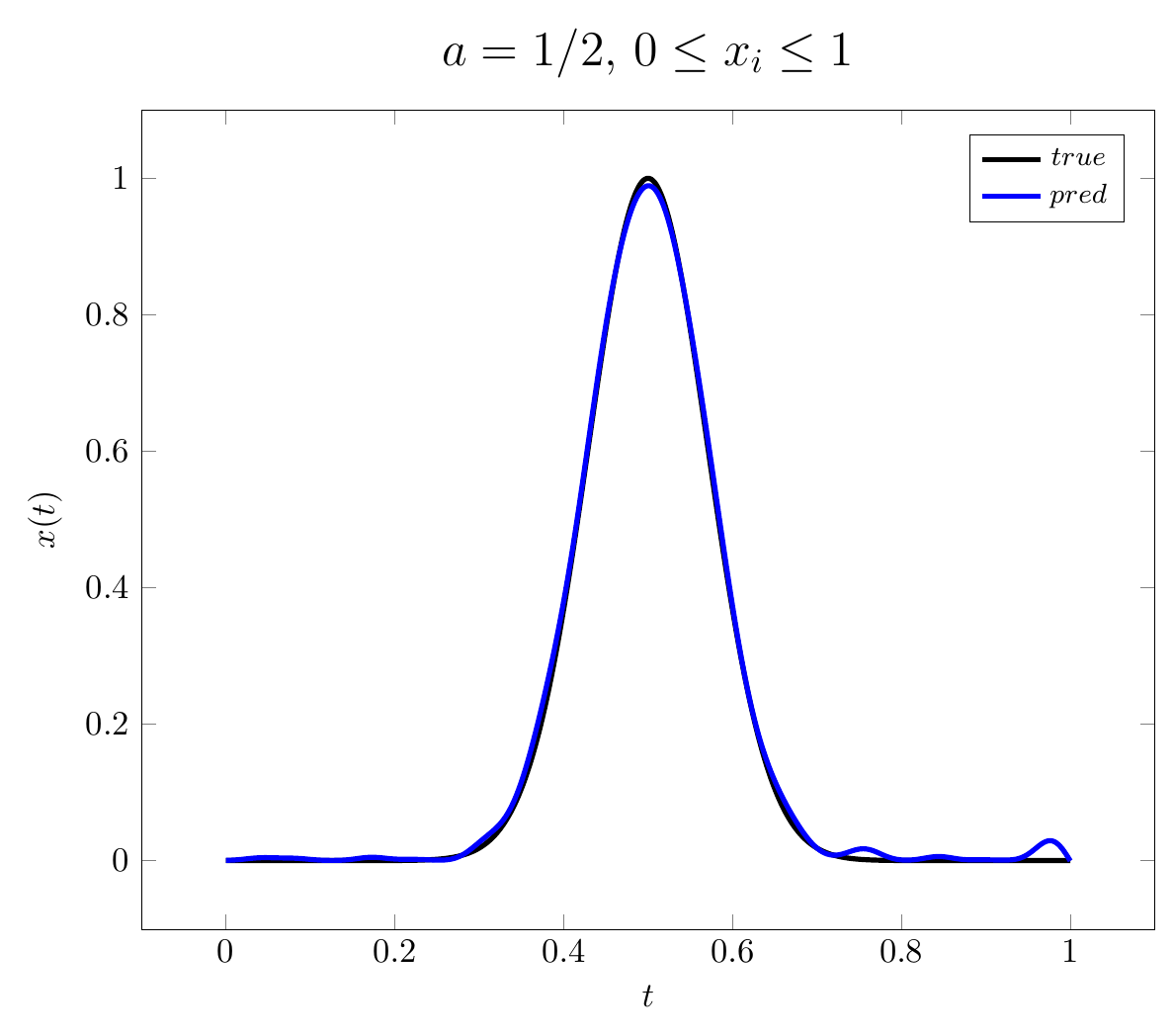}
    \includegraphics[width=6cm,scale=0.5]{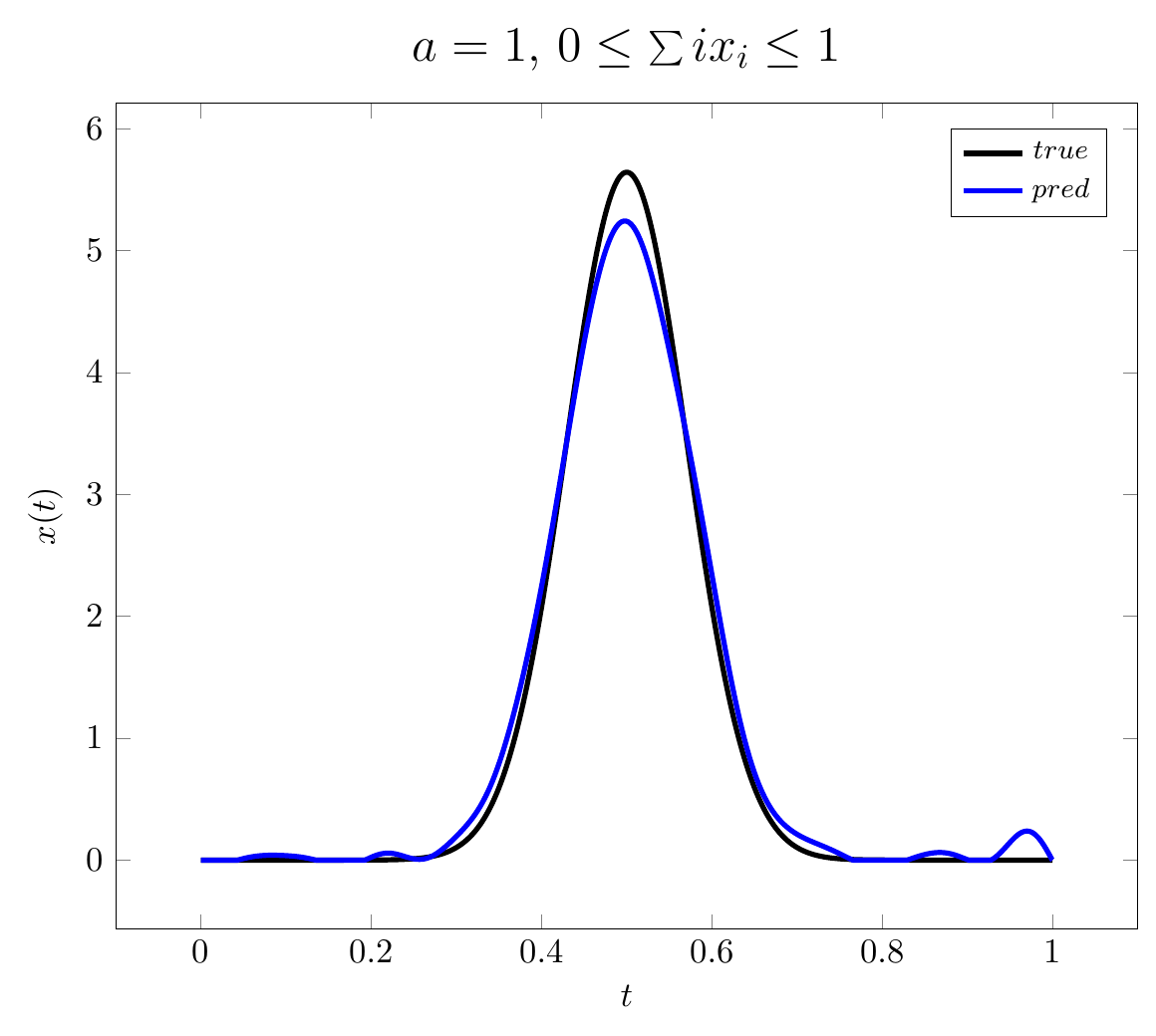}
    \includegraphics[width=6cm,scale=0.5]{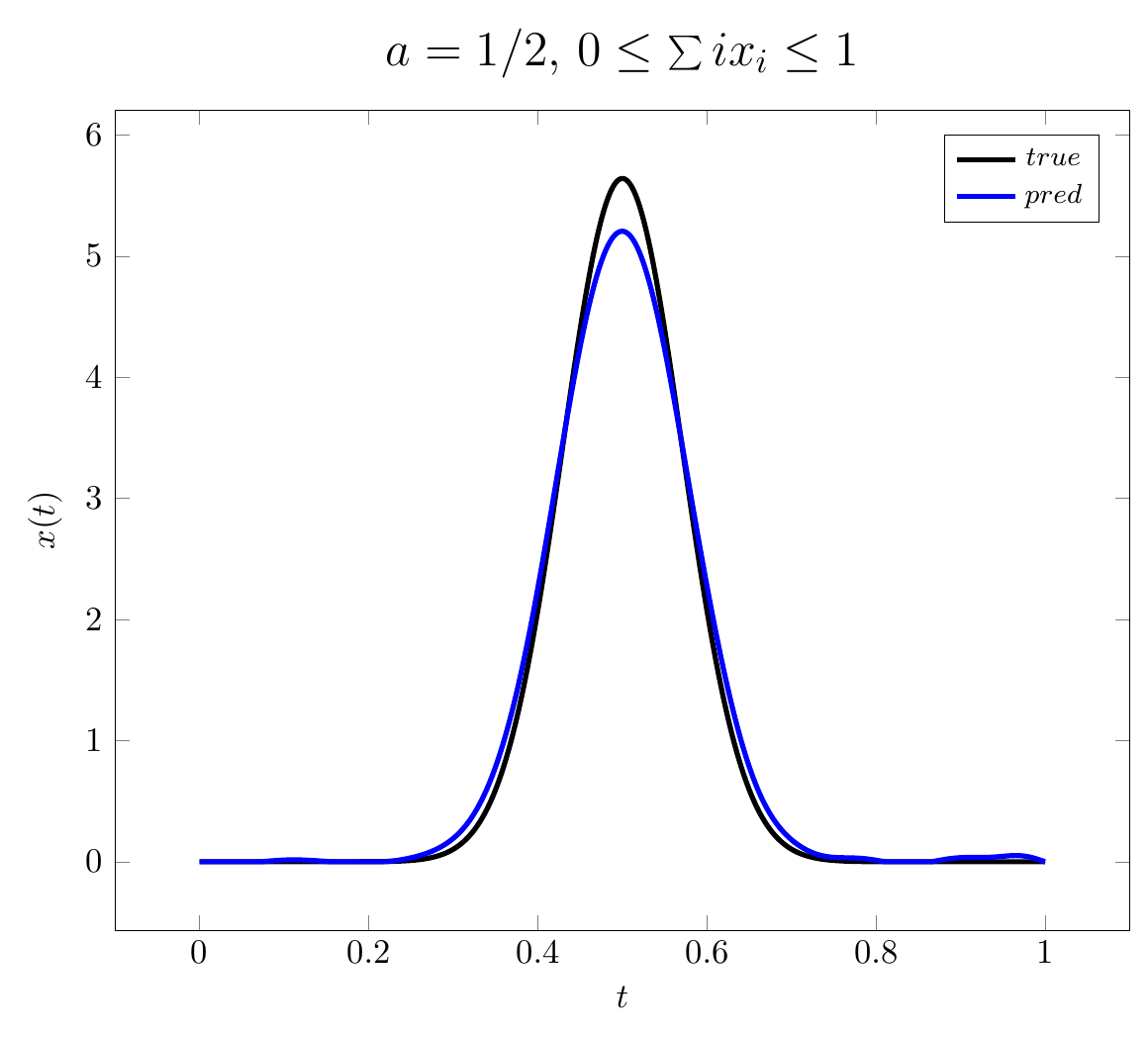}
    \caption{Output of the neural network for a Gaussian signal for various values of $a$ and various constraints. The input of the neural network convolved signal $Tx$ with an additive white noise of level $\delta = 0.05 \|Tx\|$.
    The regularization prior is based on the derivative, namely $r=1$ in~\eqref{def:varJ}, as a power of $B$ defined in~\eqref{def:D}.}
    \label{fig:gauss_pred}
\end{figure}

\begin{center}
\begin{figure}
\begin{adjustwidth}{-1.2cm}{1.2cm}  
    \begin{tabular}{ccc}
    \includegraphics[width=0.35\textwidth]{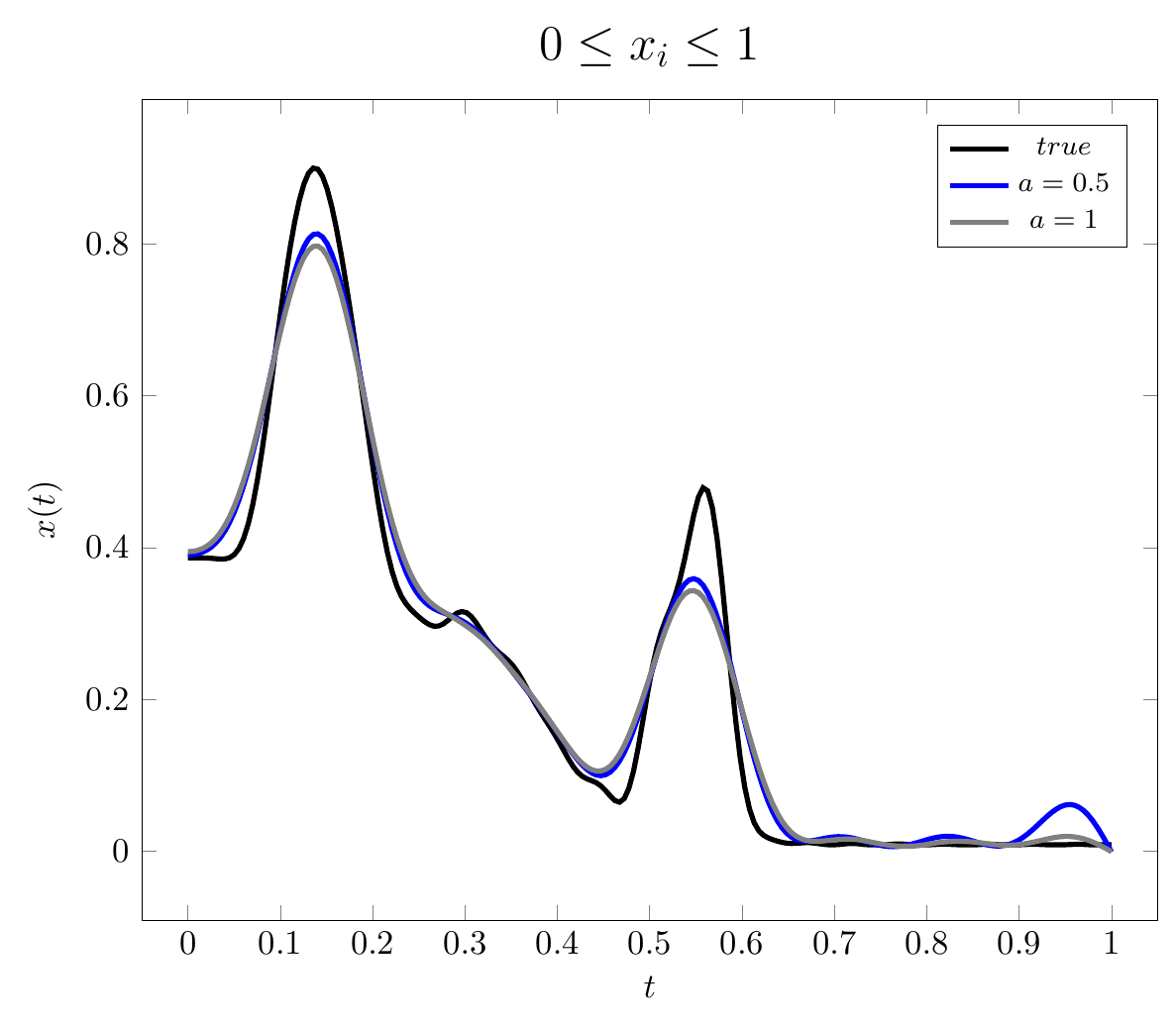}
    & \includegraphics[width=0.35\textwidth]{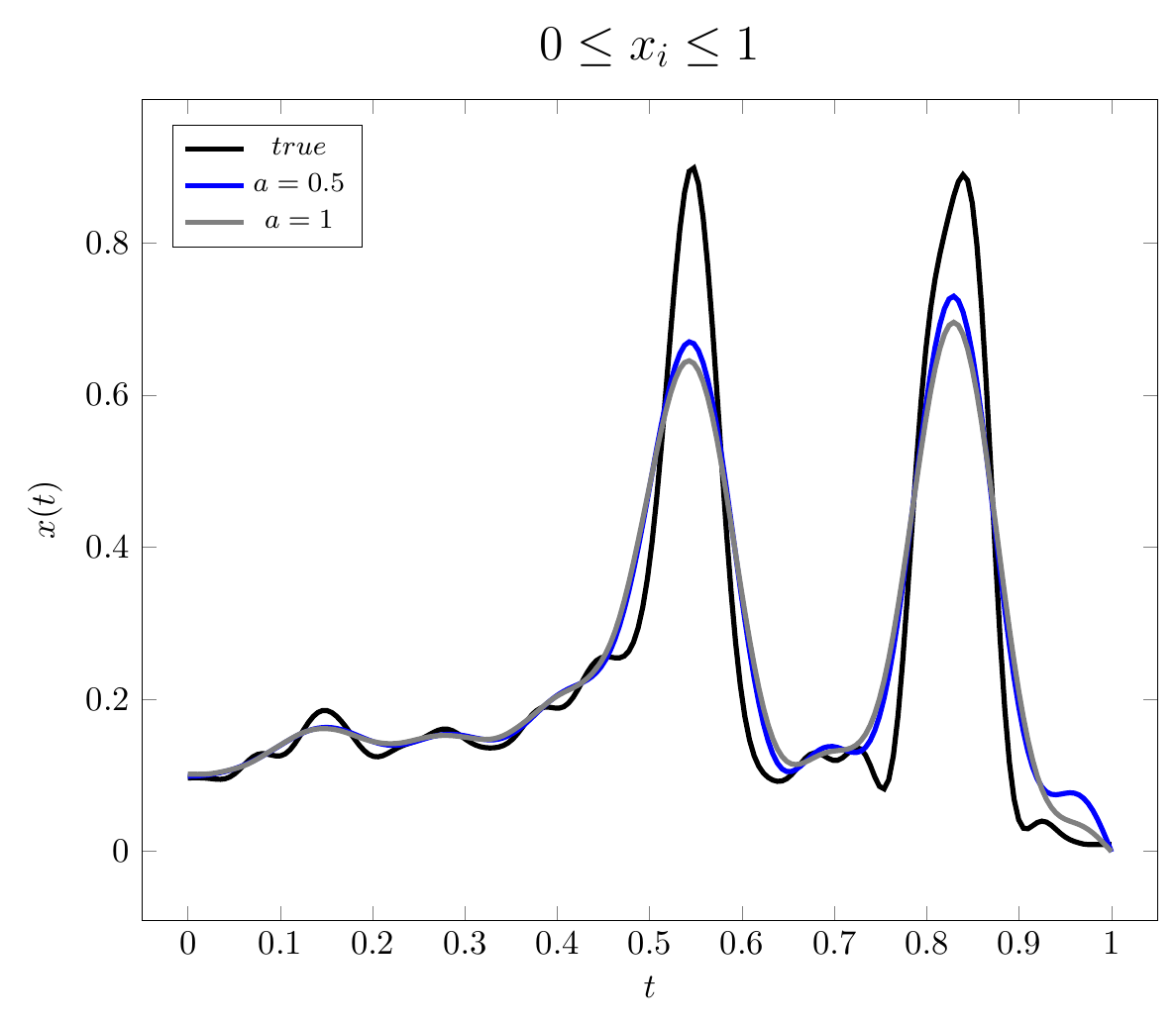}
    &\includegraphics[width=0.35\textwidth]{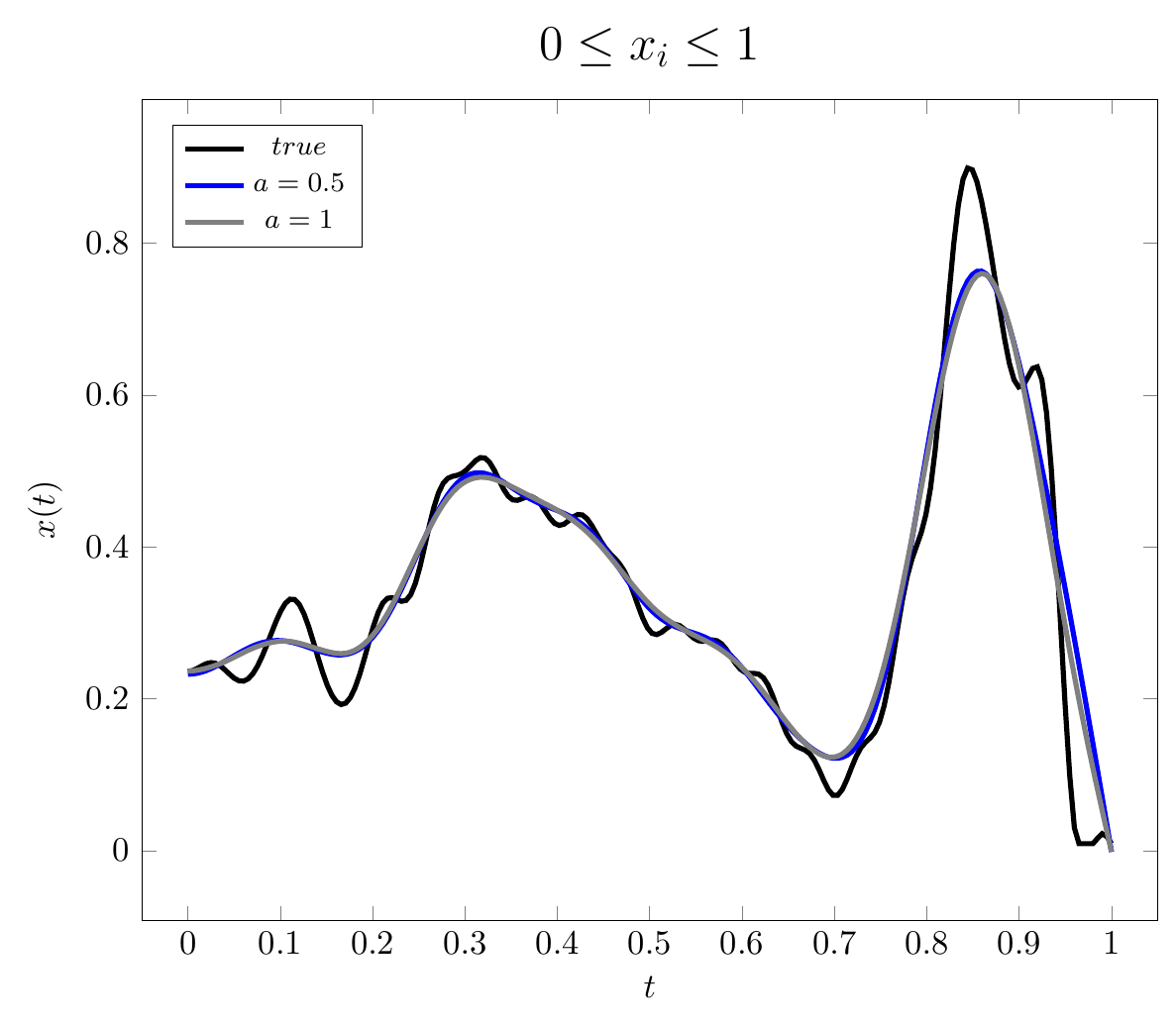} \\
    \includegraphics[width=0.35\textwidth]{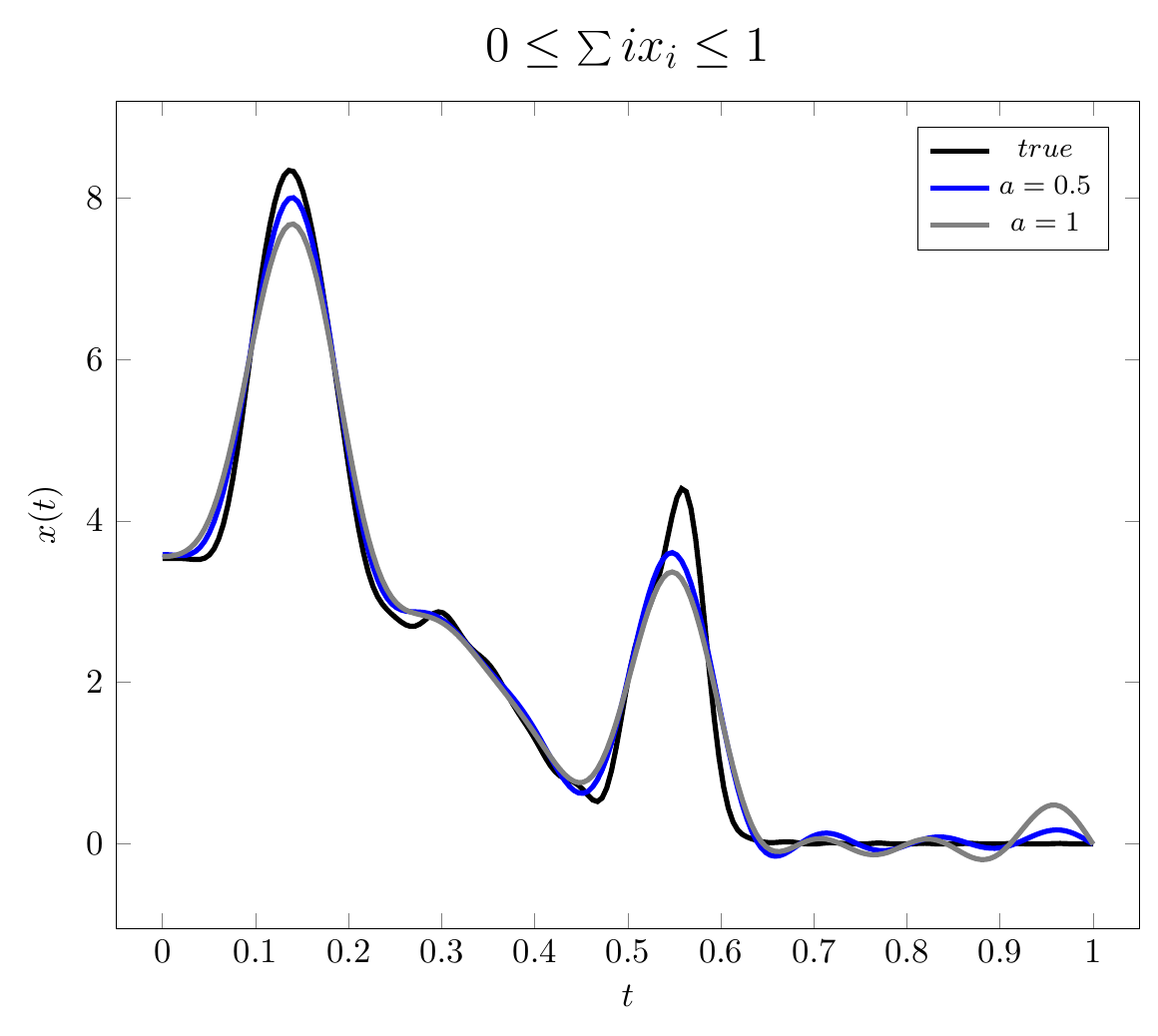}
    &\includegraphics[width=0.35\textwidth]{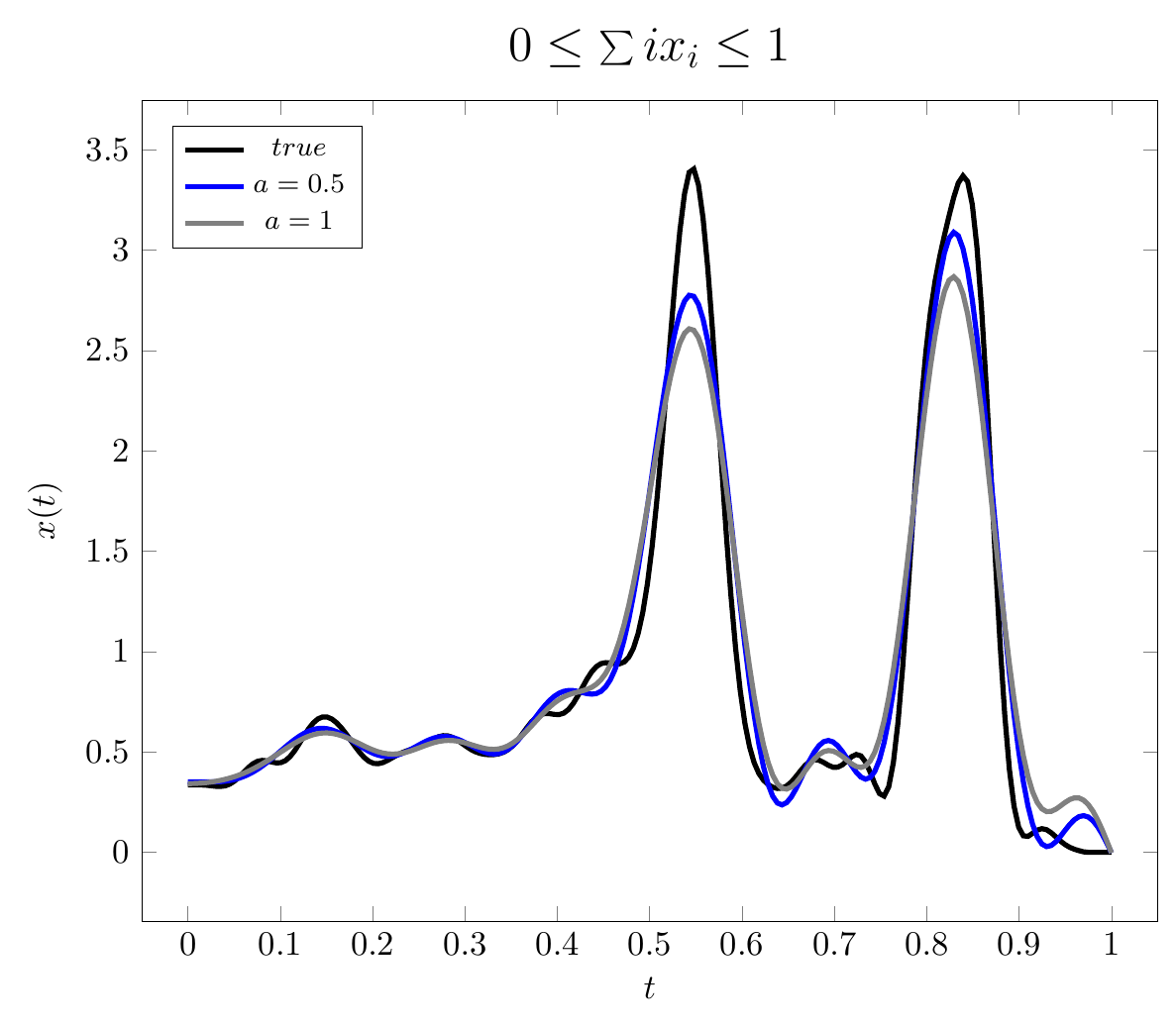}
    &\includegraphics[width=0.35\textwidth]{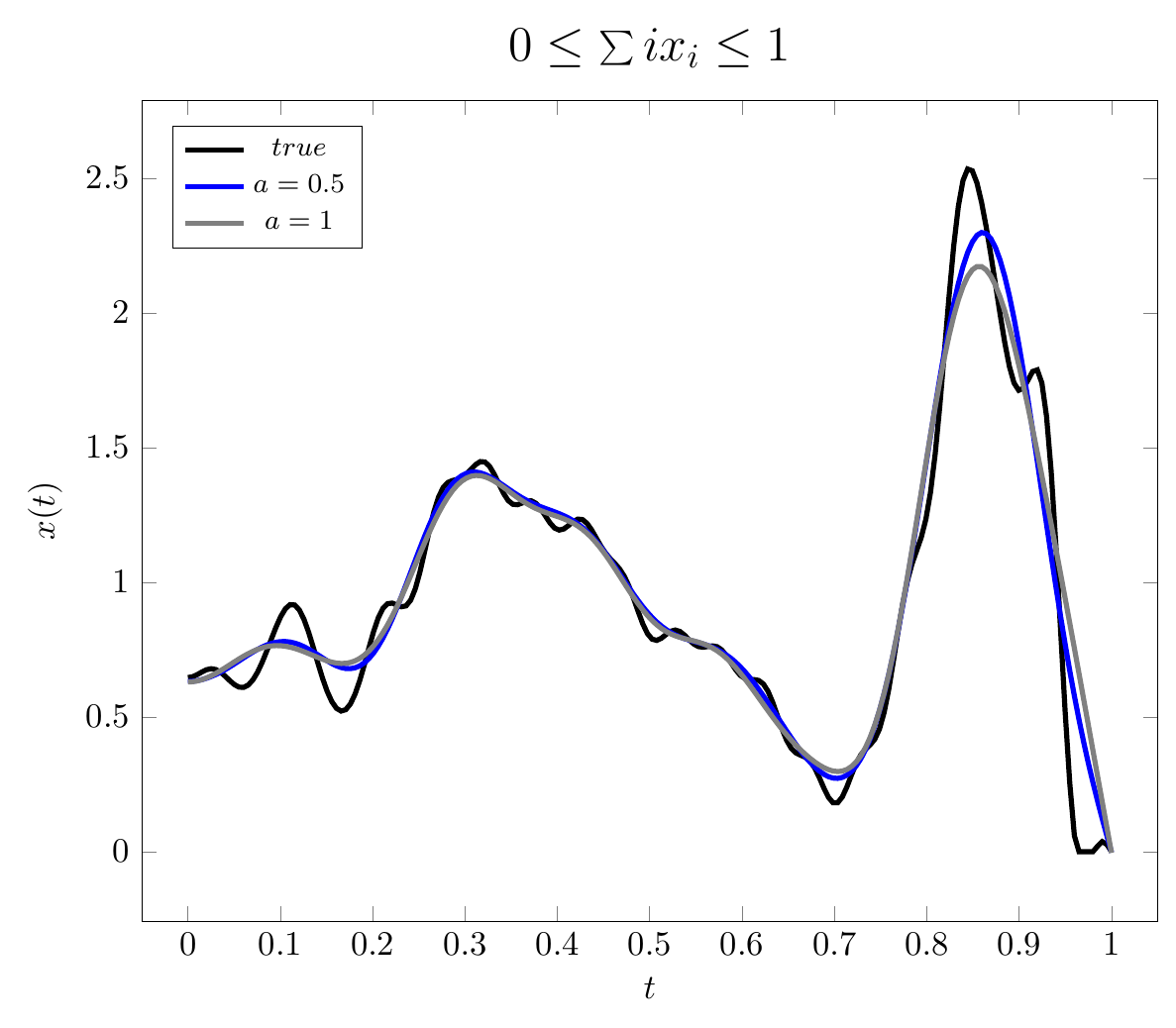}
    \end{tabular}
\end{adjustwidth}
    \caption{Example of outputs for three functions in the dataset. The measured signal is presented in Figures~\ref{fig:yydelta-a=1} and \ref{fig:yydelta-a=1/2}. The constraint $ 0 \leq \sum_i x_i \leq 1$ seems to give better result. We can also compare $a=1$ and $a=1/2$: when the order is smaller, the outputs look closer to the true signal.}
    \label{fig:diff-signal}
\end{figure}
\end{center}

\paragraph{Results}
We display in Figures~\ref{fig:gauss_pred} and~\ref{fig:diff-signal} the output of the neural network for respectively Gaussian signals and signals displayed in Figure~\ref{fig:dataset}, different values of $a$ and two possible choices for set $C$ defined either by~\eqref{constr:cube} or by~\eqref{constr:slab}.
We notice that the method performs well, and that under identical conditions, the obtained signal tends to be of lower quality when the order of the inverse problem $a$ increases. 
The reason may be numerical. If we compare this performance to a classical gradient descent algorithm, when $a$ increases the eigenvalues are greater, 
the gradient norm increases and, in order to compensate, the gradient stepsize decreases. 
Since the number of iteration is fixed, this could mean that the solution is away from the optimal solution of~\eqref{def:varJ0}.
The reason can also be theoretical: the convergence rate of the error with respect to the noise standard deviation decreases as $a$ increases according to~\cite{natterer1984}. 
\begin{figure}[h]
    \centering
    \includegraphics[width=4cm,scale=0.5]{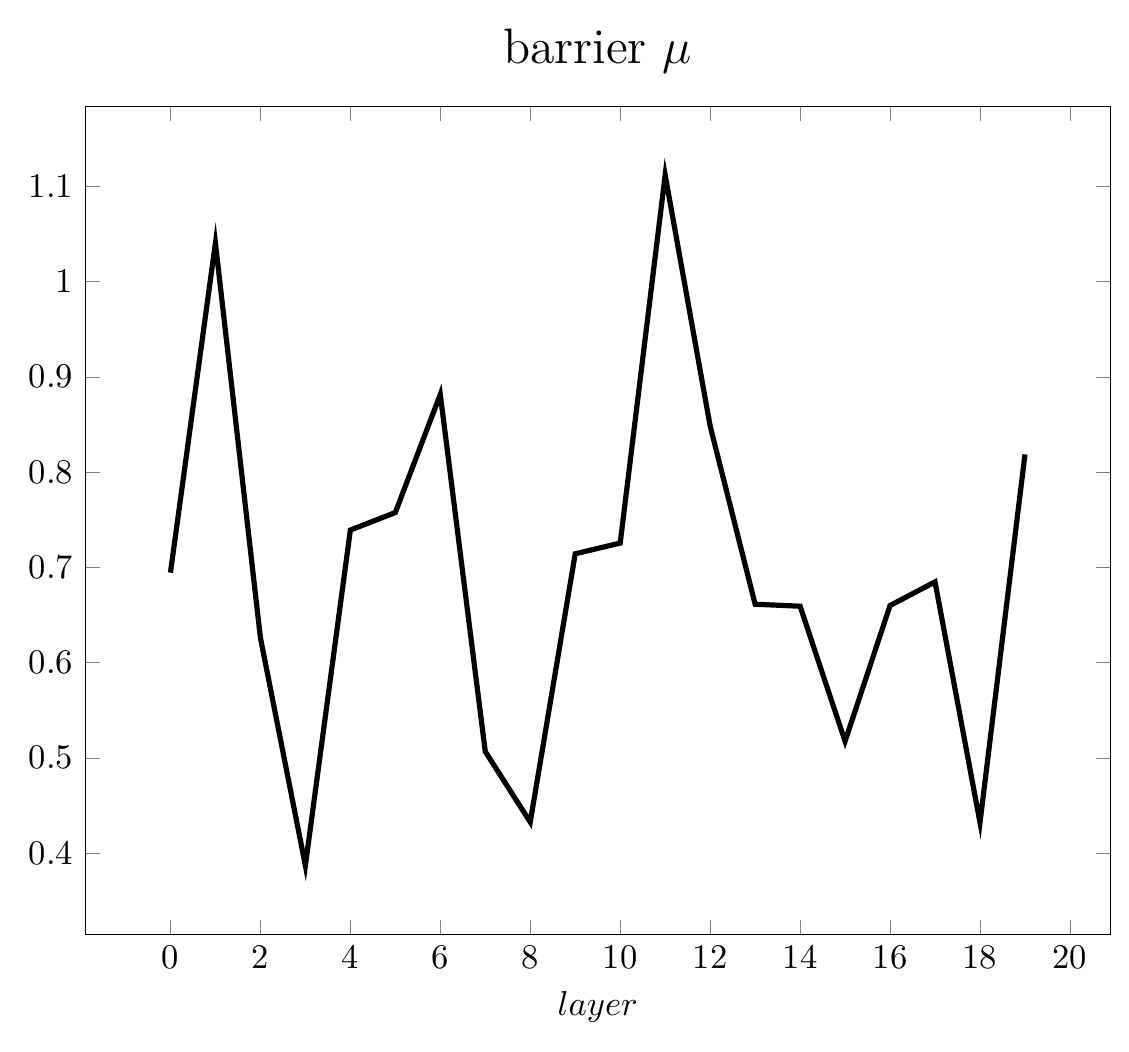}
    \includegraphics[width=4cm,scale=0.5]{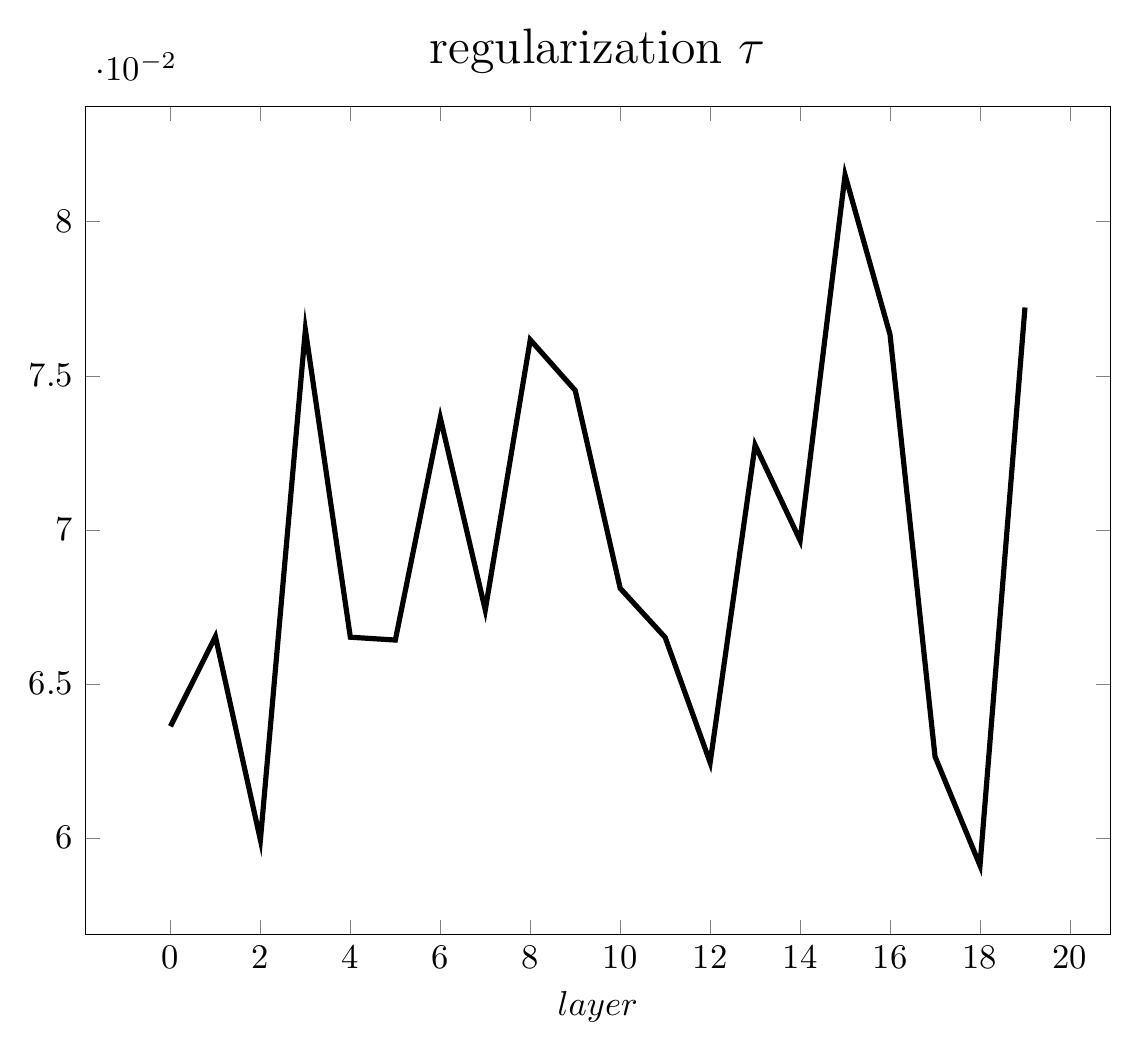}
    \includegraphics[width=4cm,scale=0.5]{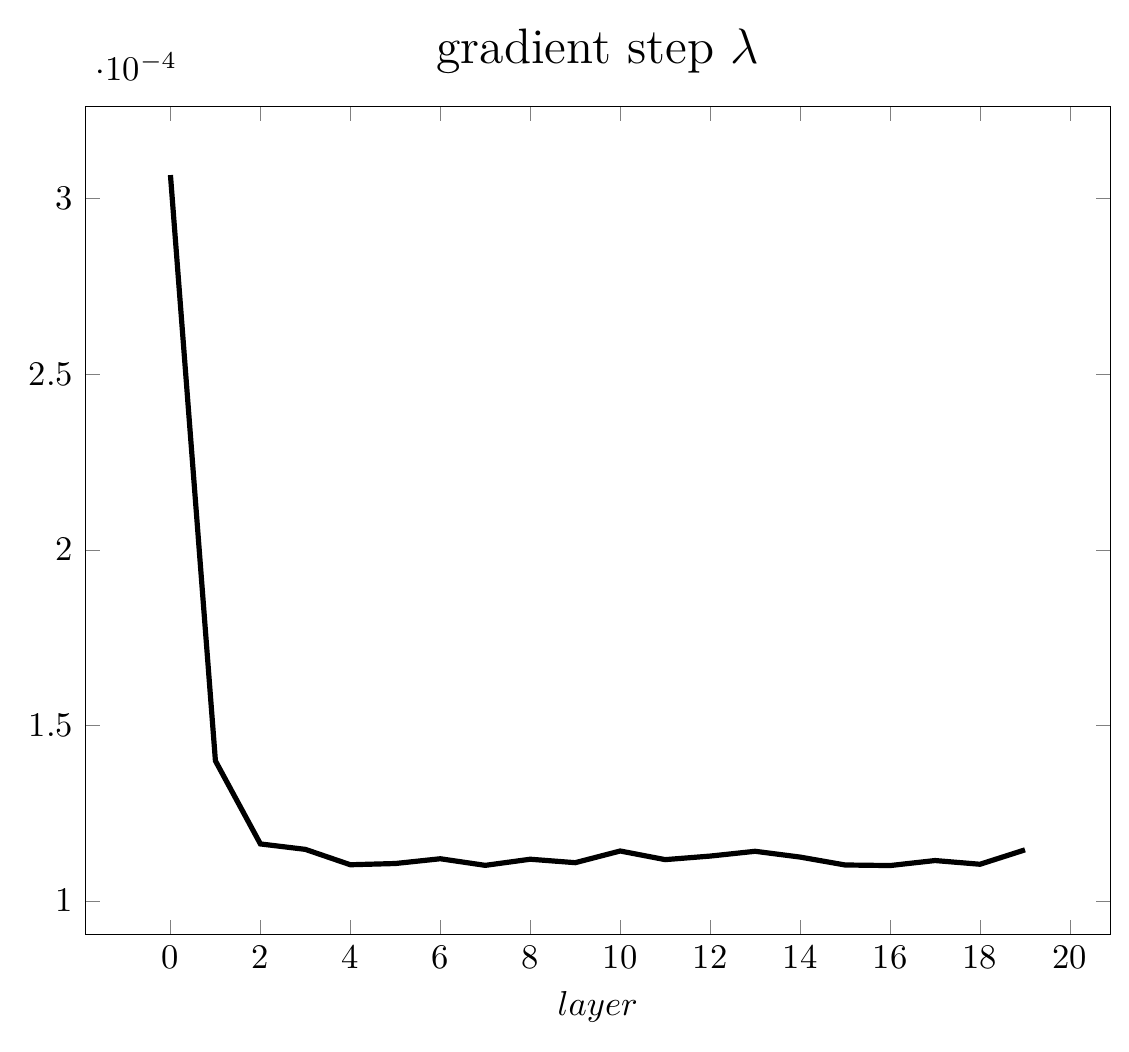}
    \caption{Hyper-parameters obtained after training for case~$(1)$ in Table~\ref{table:lip}.}
    \label{fig:hyper-parameters}
\end{figure}
\begin{table}[h]
\centering
\begin{tabular}{lcccc}
 & (1)     & (2)      &  (3)   & (4)   \\
 \hline
 & $a=1$   & $a=1/2$  & 
  $a=1$   & $a=1/2$  \\
&  $0\leq x_i \leq 1$ & $0\leq x_i \leq 1$ & $0\leq h^2 \sum_i i x_i  \leq 1$ 
& $0\leq   h^2 \sum_i i x_i \leq 1$  \\[0.2cm]
  \hline
 \eqref{e:defain-projm2}~\eqref{caseNN1} &  
$4.92\times 10^{-2}$ & $4.91\times 10^{-2}$  & 
$4.83 \times 10^{-2}$& $4.83 \times 10^{-2}$\\ 
\eqref{e:defain-projm2}~\eqref{caseNN2} &  
$2.82 \times 10^{-3}$ & $7.00 \times 10^{-3}$ & 
$8.72 \times 10^{-3}$ & $3.37 \times 10^{-3}$ \\
  \hline
\end{tabular}
\caption{Lipschitz constant of trained neural network~\eqref{def:modelNN} for various  choices of order $a$ and constraints, computed for an input $x_0 = b_0 = T^*(T x +v^\delta)=T^*\yd$. We recall that $h$ is the mesh stepsize equal to $1/N$. \label{table:lip}}
\end{table}
Figure~\ref{fig:hyper-parameters} shows the values of the hyper-parameters 
of the network after training for the constraint~\eqref{constr:cube} and the integral operator, namely $a=1$. 
We notice that the gradient step is smaller than $2/(\tau \beta_{D,K}) = 2 \beta_{T,K}/\tau = 8 \times 10^{-4} $, which is theoretically the largest gradient step leading to a convergent forward-backward algorithm.

\paragraph{Lipschitz constant estimation}
Table~\ref{table:lip} shows the Lipschitz constant obtained for the trained network under various conditions.
A first remark is that, for all the studied problems, the Lipschitz constant does not vary much neither according to the choice of $a$ nor of the constraint. 
% Knowing that these coefficients were calculated for the same input,
% therefore for the same noise level and the same regularity, 
% it is consistent that the gradient regularization and descent coefficients are similar, regardless of the subset of solution~$x$. The bias is $b_0 = T^*(T x + v^\delta)=\yd$ where $v^\delta$ is a  Gaussian additive white noise.
%
A second remark is that the obtained Lipschitz constants are lower than $1$.
The algorithm tends to constrict the solutions.
For $\delta =0$, the norm~$\|b_0\|$ is of the same order as~$\|T^*T x\|$.
For noiseless data, the neural network~\eqref{def:modelNN} would act as the inverse operator $(T^*T)^{-1}$.
Numerically, $1/\|T^*T \| = 1/\beta_{T,0} = 2.47$. 
This corresponds to the largest eigenvalue of $T^*T$. 
However, Lipschitz constants are smaller.
A possible explanation would be the following, regular function $x$ belong to a smaller vector space than noisy inputs $b_0$. 
When the regular signals are projected in the basis of eigenvectors of the compact operator $T^*T$, 
the coefficients decay extremely rapidly.
Numerically, for a given function of the dataset, the coefficient corresponding to the tenth eigenvalues is in average lower than $10^{-3}$ times the first coefficient.
However, the noisy input has non-zero coefficients over the entire spectrum.
We can therefore expect the neural network to behave roughly as the following filter of high frequency on the spectrum:
\[
\sum_{k=0}^{K-1} \frac{b_{0,k} }{\beta_{T,k} + \tau \beta_{D,k}} 
\; ,
\]
where $ b_{0,k} = \langle T^*T x , u_k\rangle = \beta_{T,k} x_k $, and $\tau$ is a regularization parameter.
Then, the Lipschitz constant is bounded by 
\[ \frac{1}{\beta_{T,0}} \, \text{max}_{ \, 0 \leq k \leq K-1} \, \frac{\beta_{T,k}}{\beta_{T,k} + \tau \beta_{D,k}} \; .
\]
As an example, for a Gaussian signal with mean value $0.5$ and standard deviation $\sigma =0.1$, and for $\tau$ of the order of $0.05$ (as in Figure~\ref{fig:hyper-parameters}), we obtain $L \sim \sum_k x_k \pi^2 /4(1 +\tau \, (2k+1)^4) \sim 0.035$.
This is the reason why we can expect the neural network Lipschitz constants to be 
roughly of the same order, namely around $10^{-2}$.

\paragraph{Comparison}
For any value of $a$, there are many techniques to reverse $T$ defined by~\eqref{def:T}.
For $a=1/2$, 
the Abel transform has been largely studied 
and three types of techniques are commonly used.

First we can mention interpolation techniques. 
Those consist in projecting the Abel  operator into a basis whose properties reflect the regularity of the solution. 
The interpolation can be performed using Chebyshev polynomial as in~\cite{piessens1973} or~\cite{pandey2014} or a Gaussian function set as in~\cite{Dribinski2002}.
They are fast, easy to implement and give good results for noiseless data.
However, as exact-inverse methods, they have an extreme sensitivity to noise, and a preprocessing may be needed.
Nevertheless, they have shown good properties for sparse data~\cite{sma2008}, or for not evenly distributed measurements (see~\cite{piessens1973}).

Secondly, Fourier transform techniques, which consist in projecting the signal in the Fourier basis are presented in~\cite{ma2008}, \cite{kalal1988abel} or similar techniques for any $a\leq1$ in~\cite{deHoop2017}. 
Those frequentist thresholding methods consist in reducing the weight of estimates on coefficients corresponding to smaller eigenvalues, for which the noise will overpower the signal.
Those techniques show computational efficiency and good noise rejection capabilities, 
but suffer from some drawbacks as they are accurate only for certain types of input data that have a sparse representation in the Fourier domain as shown in~\cite{kolhe2009} or more generally in~\cite{donoho1995}.

Thirdly, we can also mention Kalman techniques for optimum least-squares estimation applied to the inversion from noisy data.
Such techniques have been successfully applied to the Abel inverse transform in~\cite{hansen1985,Nunes1999}.
%
% The interpolation methods, for different bases, including Fourier representation, have already been the subject of studies and comparisons, on synthetic and real data as described in~\cite{pretzier1992}.
%

\begin{table}[h]
    \centering
    \begin{tabular}{lcccccc}
    \multirow{2}{*}{Noise $\delta$}
    &\multicolumn{2}{c}{Kalman}
    &\multicolumn{2}{c}{Neural Network}
    &\multicolumn{2}{c}{Fourier}\\
    \cline{2-3}
    \cline{4-5}
    \cline{6-7}
    &$a=1$&$a=1/2$
    &$a=1$&$a=1/2$
    &$a=1$&$a=1/2$
    \\
    \hline
    0.1 &0.436 & 0.383& 0.280 & 0.126
    & 0.237 &0.148
    \\
    0.05& 0.293& 0.260 &0.177 & 0.089
    & 0.186 &0.142\\
    0.01& 0.126 & 0.125 &0.095 & 0.075
    & 0.177 &0.140
    \end{tabular}
    \caption{Averaged normalized error of the output $\| x^\delta -x \| /\|x\|$ obtained for different noise standard deviation values $\delta$ and different types of signals.
    The error for the Kalman Filter are of the order of $\sqrt{\delta} $ according to theory, but in practice finding the parameters allowing to reach such a precision is difficult. We compare the result with the neural network with the constraint defined by~\eqref{constr:cube}.}
    \label{tab:compare}
\end{table}

Two methods have been implemented, Fourier and Kalman, 
and tested over the same dataset of 50 signals $x$,
different from the training set of the neural network. 
The averaged error on the outputs over the dataset for each methods are displayed in Table~\ref{tab:compare}.
The neural network~\eqref{def:modelNN} compares favorably with other techniques for a relative solution error in $L^2$. 
The Fourier method is more accurate as the noise level increases, since it works by filtering high frequencies. 
The Kalman method returns less regular solutions, and calibration of the regularization parameter raises an additional experimental difficulty.

\begin{figure}
    \centering
    \begin{adjustwidth}{-1.2cm}{1.2cm} 
    \begin{tabular}{ccc}
    \includegraphics[width=5cm,scale=0.5]{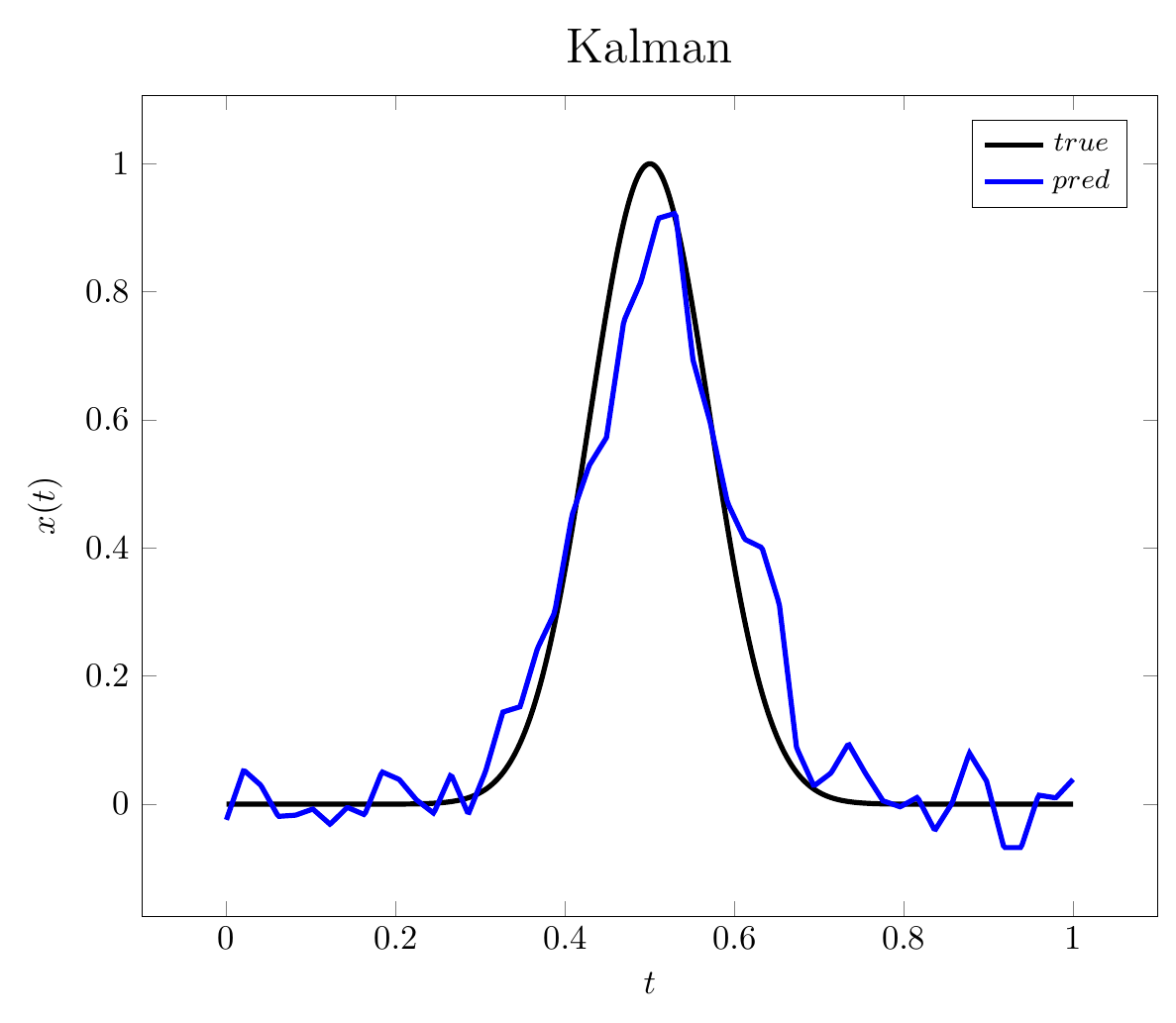}
    &\includegraphics[width=5cm,scale=0.5]{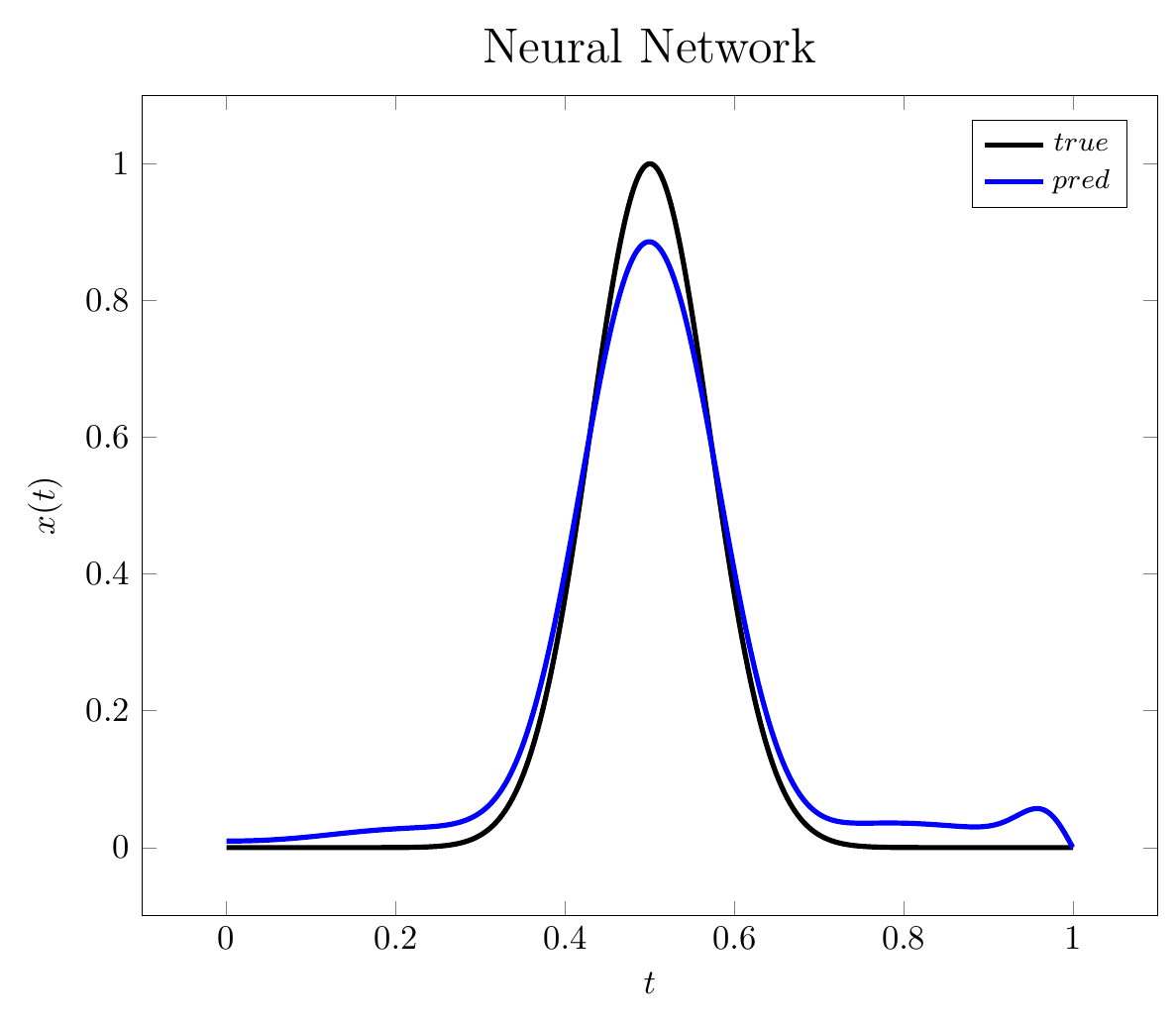}
    &\includegraphics[width=5cm,scale=0.5]{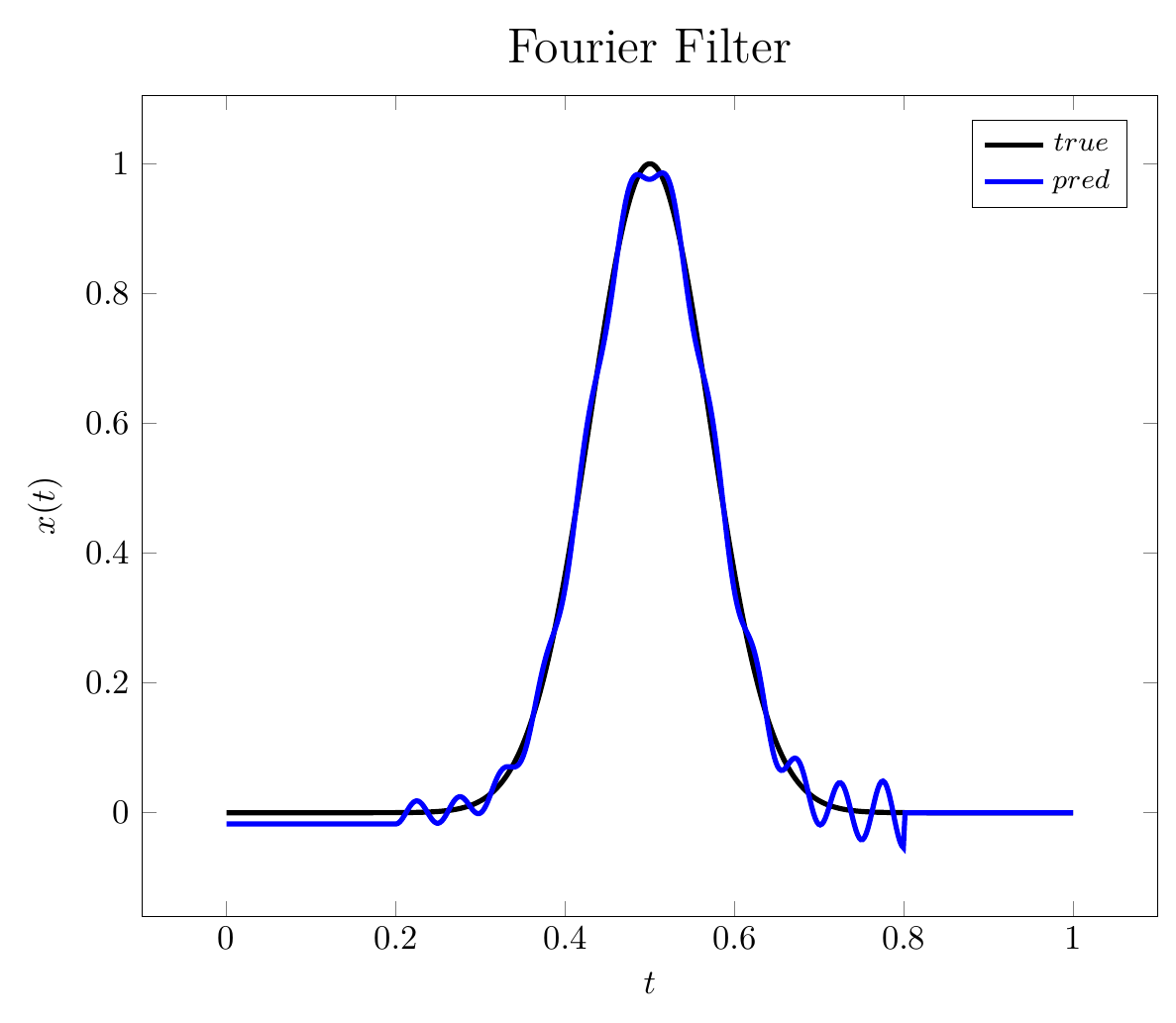}
    \end{tabular}
    \end{adjustwidth}
    \caption{Output of the neural network for a Gaussian signal when $a=1$, the white noise level is $\delta = 0.05 \|Tx\|$, using three differents techniques : Kalman, neural network with constraint~\eqref{constr:cube} and Fourier filtering.}
    \label{fig:gauss_pred2}
\end{figure}
%%%%%%%%%%%%%%%%%%%%%%%%%%%%%%%%%%%%%%%%%%%%%%%%%%%%%%%%%%%%%%%%%%%%%%%%%%%%%%%%%%%%%%%%%%%%%%%%%%%%%%%%%%%%%%%%%%%%%%%%%%%%%%%%%
% CONCLUSION
\section{Conclusion}
In the continuity of the work of~\cite{Corbineau2020}, the present paper proposes to unroll an algorithm obtained from a variational formulation of 1D integral inverse problems. 
This approach is versatile, since it allows to invert a broad family of integral or convolution operators, 
and it delivers a solution taking into account physical constraints of the problem.
Indeed, with some existing methods, it may be difficult to enforce constraints, such as complying with some bounds or belonging to a given subspace. In particular, in many practical scenarios, solutions that do not fulfill basic constraints such as positivity 
may appear as irrelevant in terms of physical interpretation.
The numerical solutions obtained for the case of the Abel operator indicate that the approach is easy to implement and computationally very attractive, since the training takes only a couple of hours and testing or prediction takes a few seconds on a regular CPU. \\
We additionally performed a theoretical analysis of robustness with respect to the observed data, which ensures the reliability of the proposed inverse method.
In future work, more sophisticated neural network structures could be considered or additional parameters (such as the leakage factors we introduced in our analysis) could be learnt. Also, training sets which would better suited to specific applications could be employed within our framework. We think also that the $\alpha$-averaged properties that we established pave the way for building recurrent networks in the spirit of~\cite{Combettes2019}.

%%%%%%%%%%%%%%%%%%%%%%%%%%%%%%%%%%%%%%%%%%%%%%%%%%%%%%%%%%%%%%%%%%%%%%%%%%%%%%%%%%%%%%%%%%%%%%%%%%%%%%%%%%%%%%%%%%%%%%%%%%%%%%%%%

\bibliographystyle{unsrt}%Used BibTeX style is unsrt
\bibliography{biblio}
\end{document}